# COMBINATORIAL LÉVY PROCESSES

HARRY CRANE

ABSTRACT. Combinatorial Lévy processes evolve on general state spaces of countable combinatorial structures. In this setting, the usual Lévy process properties of stationary, independent increments are defined in an unconventional way in terms of the symmetric difference operation on sets. In discrete time, the description of combinatorial Lévy processes gives rise to the notion of combinatorial random walks. These processes behave differently than random walks and Lévy processes on other state spaces. Standard examples include processes on sets, graphs, and $n$-ary relations, but the framework permits far more general possibilities. The main theorems characterize both finite and infinite state space combinatorial Lévy processes by a unique $\sigma$-finite measure. Under the additional assumption of exchangeability, we obtain a more explicit characterization by which every exchangeable combinatorial Lévy process corresponds to a Poisson point process on the same state space. Associated behavior of the projection into a space of limiting objects reflects certain structural features of the covering process.

## 1. INTRODUCTION

A Lévy process $(X_t, t \geq 0)$ on $\mathbb{R}^d$ is a random map $t \mapsto X_t$ with stationary, independent increments and càdlàg sample paths with respect to the Euclidean topology. Lévy processes comprise a large class of tractable models with applications in finance, neuroscience, climate modeling, etc. The Lévy–Itô–Khintchine theorem decomposes their rich structure into an independent Brownian motion with drift, a compound Poisson process, and a pure jump martingale. Bertoin [7] surveys the properties of $\mathbb{R}$-valued Lévy processes, which specialize those of Lévy processes in topological groups. In an arbitrary topological group $\mathcal{X}$, the Lévy process assumptions are defined with respect to the group action, with the left, respectively right, increment between $x, x' \in \mathcal{X}$ defined as the unique $y \in \mathcal{X}$ such that $x = yx'$, respectively $x = x'y$. Liao [33] gives a general introduction to Lévy processes in topological groups with special treatment of the Lie group case, which garners interest for its relation to certain types of stochastic flows.

In both the real-valued and Lie group settings, many nice properties result from the interplay between the increments assumptions and the topology of the underlying state space. In Euclidean space, the Lévy–Itô–Khintchine representation is tied to its predecessor, the Lévy–Khintchine theorem for infinitely divisible distributions. In a Lie group, the smoothness of the associated Lie algebra plays a key role.

Afield of Lévy processes, combinatorial stochastic processes evolve on discrete state spaces, with a focus on the theory of exchangeable random partitions [16, 23, 31], coalescent and fragmentation processes [6, 8, 32, 36], connections to stable subordinators, Brownian bridges, and Lévy processes [35, 38], tree- [1, 2, 3, 11, 22, 37] and graph-valued [13, 14, 15] processes.

Below I introduce *combinatorial Lévy processes* as a family of models that may be suitable for dynamic structures that arise in streaming data collection, complex networks, and other applications. Combinatorial Lévy processes evolve on discrete spaces of labeled combinatorial objects with the following special cases as an illustration.

- *Set-valued processes*: On the space of subsets of $\mathbb{N} := \{1, 2, \ldots\}$, a combinatorial Lévy process evolves by rearranging elements. For example, each element $i = 1, 2, \ldots$ might enter and leave the set at alternating times of independent rate-1 Poisson processes. More generally, each element $i = 1, 2, \ldots$ can enter and leave at the alternating times of independent rate-$i$ Poisson processes, making the behavior of different elements inhomogeneous. These dynamics imitate those of some previously studied partition-valued processes, for example, [9, 12, 35]. Forty years ago, Harris [25, 26] studied set-valued processes under entirely different assumptions.
- *Graph-valued processes*: Dynamic networks arise in a range of modern applications involving time-varying interactions in a population, for example, [24, 27, 29, 39]. Our main discussion describes the possibilities and limitations of combinatorial Lévy process models for dynamic networks. The Lévy process assumptions, if appropriate, make these processes particularly applicable in statistical applications, as estimating the increments distribution is a straightforward computational exercise. We discuss these statistical aspects further in Section 8.1.
- *Networks with community structure*: The most interesting context for the theory of combinatorial Lévy processes is in joint modeling of composite structures, such as dynamic networks with underlying communities of vertices. In this case, it is natural and straightforward to combine the above two processes on sets and graphs into one that models the joint evolution of a network and a community of its vertices. Extensions to collections of $k$ different communities and [$l$] different networks and possibly higher order interactions also fall within the scope of combinatorial Lévy processes. These processes incorporate temporal variation into the widespread statistical literature on community detection in networks, almost all of which is confined to the case in which the underlying network is static. This example, therefore, provide a concrete context within which to illustrate the more general theory of combinatorial Lévy processes developed below.

1.1. **Relationship to other literature.** The theory of combinatorial Lévy processes introduced here follows a line of research on combinatorial stochastic process models for various applications, see [9, 35], but under different assumptions and in a more general context.

Though we often illustrate the main concepts with examples of set- and graph-valued processes, these are not the main intended applications of combinatorial Lévy processes, and prior works in either of these settings do not consider any analog to the Lévy process assumptions introduced below. Applications such as the last one above, in which the joint evolution of a community structure and a network is modeled, more fully capture the spirit of combinatorial Lévy processes, as they permit dependence in the evolution of disjoint components. In this regard, certain combinatorial Lévy processes may serve as suitable null models for epidemic spread on dynamic networks, extending prior work in [24] and serving as a tractable setting in which to address some questions posed in the concluding remarks of [21].

We also note that the special cases of sets and graphs, though helpful as examples, do not exhibit all of the subtleties at play in the more general theorems. In particular, the general results for combinatorial Lévy processes presented here are disjoint from other recent work that specializes to the case of exchangeable graph-valued Feller processes [15]. Section 5.2 shows the main qualitative similarity between the graph-valued processes of [15] and Lévy processes evolving on the space of undirected graphs. Heuristically, the decomposition of jump rates in Section 5.2 corresponds roughly to the decomposition in [15, Sec. 3.2.1], but otherwise the main theorems in [15] differ in several important respects from our main theorems for combinatorial Lévy processes.

Most notably from the standpoint of application, the starting points of [15] and our study of combinatorial Lévy processes are disjoint. The techniques in [15] rely critically on the assumptions that the state space consists of all countable graphs and that the processes are exchangeable and exhibit a property of Markovian consistency under subsampling. Except those theorems below that explicitly assume infinite exchangeability, our results rely on none of these prior assumptions. We instead assume a process, which may be exchangeable or not, on a space of arbitrary combinatorial objects, which may be labeled in a finite or countable set, whose sample paths exhibit the properties of stationary and independent increments defined with respect to a group action on the corresponding space. Because of these differences, our main theorems hinge on different observations and techniques.

Specifically, Theorems 4.5 and 4.6 uniquely characterize the behavior of combinatorial Lévy processes without the requirement of exchangeability. Moreover, the characteristic measure of combinatorial Lévy processes in Theorems 4.5, 4.6, 4.15, and 4.18 is uniquely determined by the process **X**, whereas the characterization of exchangeable Feller processes in [15, Theorem 3] is not unique. Perhaps most significantly, the characteristic measures in [15, Theorems 3 and 4] are defined on an abstract space of *rewiring maps* which differs from the state space of graphs on which the process is defined. Our main theorems for combinatorial Lévy processes on $\mathcal{L}_\mathbb{N}$, on the other hand, characterize the behavior by a unique measure on the same space $\mathcal{L}_\mathbb{N}$. On the space of graphs, therefore, the characteristic measure for combinatorial Lévy processes is also defined as a measure on the space of graphs, which differs from [15] even in this special case. Finally, as alluded above and discussed further in Section 8.1, the representation of combinatorial Lévy processes given below offers a straightforward way to estimate the transition measure from data. The characterization in [15], on the other hand, is given in terms of a measure on *rewiring maps*, which is far more difficult to estimate from data and, therefore, less applicable than combinatorial Lévy processes. The same comparison holds for the more abstract discussion of Markov processes on a Fraïssé space studied more recently in [17].

1.2. **Outline.** In Section 2, we summarize the main theorems in the case of set-valued Lévy processes. In Section 3, we lay down key definitions, notation, and observations. In Section 4, we formally summarize the main theorems in the language of Section 3. In Section 5, we demonstrate our main theorems with concrete examples that are relevant to specific applications. In Section 6, we prove a key theorem about $\sigma$-finite measures on combinatorial spaces, from which we readily deduce the Lévy–Itô–Khintchine representation for exchangeable combinatorial Lévy processes. In Section 7, we prove our main theorems. In Section 8, we make concluding remarks.

2. Exposition: set-valued processes

**Remark 2.1** (Notation). *We discuss both discrete and continuous time processes. When speaking generally, we index time by $t \in [0, \infty)$. When speaking specifically about discrete time processes, we index time by $m \in \mathbb{Z}_+ := \{0, 1, \ldots\}$ and write $\mathbf{X} = (X_m, m \geq 0)$ to denote a discrete time process.*

The concept of combinatorial increments captures structural differences between combinatorial objects. To fix ideas, we first assume that $\mathbf{X} = (X_t, t \geq 0)$ evolves on the space of subsets of a base set $S \subseteq \mathbb{N}$, denoted $2^S$.

2.1. **Increments and topology.** Every $A \subseteq \mathbb{N}$ determines a map $2^{\mathbb{N}} \to 2^{\mathbb{N}}$ by $A' \mapsto A \triangle A'$, where

(1) $$A \triangle A' := (A \cap A'^c) \cup (A^c \cap A')$$

is the *symmetric difference* operation and $A^c := \mathbb{N} \setminus A$ denotes the *complement* of $A$ (relative to $\mathbb{N}$). Under this operation, the empty set $\emptyset := \{\}$ acts as the identity and each $A \subseteq \mathbb{N}$ is its own inverse, that is, $A \triangle A = \emptyset$ for all $A \subseteq \mathbb{N}$. We equip $2^{\mathbb{N}}$ with the product discrete topology induced by

(2) $$d(A, A') := 1/(1 + \sup\{n \in \mathbb{N} : A \cap [n] = A' \cap [n]\}), \quad A, A' \subseteq \mathbb{N},$$

with the convention $1/\infty = 0$, where $[n] := \{1, \ldots, n\}$. As we do for arbitrary combinatorial spaces below, we equip $2^{\mathbb{N}}$ with its Borel $\sigma$-field under (2) and, if $S \subseteq \mathbb{N}$ is infinite, we equip $2^S$ with the trace of the Borel $\sigma$-field on $2^{\mathbb{N}}$. If $S \subset \mathbb{N}$ is finite, then we equip $2^S$ with the discrete metric and the power set $\sigma$-field $2^{2^S}$. Any subset $\mathcal{B} \subseteq 2^S$ that we discuss below is implicitly assumed to be measurable with respect to the corresponding $\sigma$-field.

In the following definition, $T$ stands for either discrete time ($T = \mathbb{Z}_+$) or continuous time ($T = [0, \infty)$). The definition holds in either case, the only difference being that càdlàg paths are automatic in discrete time.

**Definition 2.2** (Combinatorial Lévy process on $2^S$). *We call $\mathbf{X} = (X_t, t \in T)$ a combinatorial Lévy process on $2^S$ if it has*
- $X_0 = \emptyset$,
- *stationary increments, that is, $X_{t+s} \triangle X_s =_{\mathcal{D}} X_t$ for all $s, t \geq 0$, where $=_{\mathcal{D}}$ denotes equality in distribution,*
- *independent increments, that is, $X_{t_1} \triangle X_{t_0}, \ldots, X_{t_k} \triangle X_{t_{k-1}}$ are independent for all $0 \leq t_0 \leq t_1 \leq \cdots \leq t_k < \infty$ in $T$, and*
- *càdlàg sample paths, that is, $t \mapsto X_t$ is right continuous and has left limits under the induced topology on $2^S$.*

In discrete time, we interpret a combinatorial Lévy process on $2^S$ as a set-valued random walk.

**Definition 2.3** (Set-valued random walk). *A random walk on $2^S$ with increment distribution $\mu$ on $\mathcal{L}_{\mathbb{N}}$ and initial state $X_0$ is a discrete time process $\mathbf{X} = (X_m, m \geq 0)$ with $X_{m+1} =_{\mathcal{D}} X_m \triangle \Delta_{m+1}$ for every $m \geq 0$, where $\Delta_1, \Delta_2, \ldots$ are independent and identically distributed (i.i.d.) according to $\mu$.*

**Theorem 2.4.** *For any $S \subseteq \mathbb{N}$, let $\mathbf{X} = (X_m, m \geq 0)$ be a discrete time combinatorial Lévy process on $2^S$. Then there exists a unique probability measure $\mu$ on $2^S$ such that $\mathbf{X}$ is distributed as a random walk with initial state $\emptyset$ and increment distribution $\mu$.*

The proof of Theorem 2.4 is straightforward even for general combinatorial Lévy processes—see Theorem 4.5—but we explicitly prove the set-valued case to aid the more general discussion below.

*Proof.* The stationary and independent increments assumptions imply that $\mathbf{X} = (X_m, m \geq 0)$ is determined by its initial state $X_0 = \emptyset$ and an independent, identically distributed (i.i.d.) sequence $\mathbf{\Delta} = (\Delta_m, m \geq 1)$ of subsets, where

$$\Delta_m = X_m \triangle X_{m-1}, \quad m \geq 1.$$

For each $m \geq 1$, $\Delta_m$ contains all elements whose status changes between times $m-1$ and $m$; thus, the transition law of $\mathbf{X}$ is governed by a unique probability measure $\mu$ on $2^S$, which acts as the increments measure for the random walk started at $\emptyset$. $\square$

In continuous time, we must distinguish between processes on structures labeled by finite and infinite sets. In particular, a process $\mathbf{X} = (X_t, t \geq 0)$ on $2^{\mathbb{N}}$ can experience infinitely many jumps in bounded time intervals, but càdlàg sample paths constrain each induced finite state space process $\mathbf{X}^{[n]} := (X_t \cap [n], t \geq 0)$ to jump only finitely often in bounded intervals. These competing notions harness the behavior of $\mathbf{X}$ in a special way. The same behavior persists for $\mathbf{X}$ on $2^S$ for any infinite set $S \subseteq \mathbb{N}$, but we lose no generality in assuming $S = \mathbb{N}$ in this case.

First, we observe that any combinatorial Lévy process has the Feller property (Corollary 4.10) and, thus, its evolution is determined by the infinitesimal jump rates

$$\mu(d\Delta) = \lim_{t \downarrow 0} \frac{1}{t} \mathbb{P}\{X_t \in d\Delta\}, \quad \Delta \in 2^{\mathbb{N}} \setminus \{\emptyset\}.$$

Since $\Delta = \emptyset$ corresponds to no jump, we may tacitly assume $\mu(\{\emptyset\}) = 0$. To ensure that each $\mathbf{X}^{[n]}$ jumps only finitely often, $\mu$ must also satisfy

$$\mu(\{A \in 2^{\mathbb{N}} : A \cap [n] \neq \emptyset\}) < \infty \quad \text{for all } n \in \mathbb{N}.$$

Since the behavior of $\mathbf{X}$ is determined by the infinitesimal jump rates $\lim_{t \downarrow 0} t^{-1} \mathbb{P}\{X_t \in d\Delta\}$, $\mathbf{X}$ can be described by a unique measure $\mu$ on $2^{\mathbb{N}}$ that satisfies (3) below. We summarize these observations with the following general statement for arbitrary set-valued Lévy processes on finite and infinite state spaces.

**Theorem 2.5.** *Let $\mathbf{X} = (X_t, t \geq 0)$ be a continuous time combinatorial Lévy process on $2^S$. Then there is a unique measure $\mu$ on $2^S$ for which*

(3) $$\mu(\{\emptyset\}) = 0 \quad \text{and} \quad \mu(\{A \in 2^S : A \cap [n] \neq \emptyset\}) < \infty \quad \text{for all } n \in \mathbb{N}$$

*such that the infinitesimal jump rates of $\mathbf{X}$ satisfy*

$$\lim_{t \downarrow 0} \frac{1}{t} \mathbb{P}\{X_t \in d\Delta\} = \mu(d\Delta), \quad \Delta \in 2^S \setminus \{\emptyset\}.$$

**Remark 2.6.** *Note that if $S \subset \mathbb{N}$ is finite, then the righthand side of (3) implies that $\mu$ is a finite measure.*

From any $\mu$ satisfying (3), we construct the *$\mu$-canonical Lévy process* $\mathbf{X}^*_\mu = (X^*_t, t \geq 0)$ from a Poisson point process $\mathbf{\Delta}^* = \{(t, \Delta_t)\} \subseteq [0, \infty) \times 2^S$ with intensity measure $dt \otimes \mu$, where $dt$ denotes Lebesgue measure on $[0, \infty)$. The atoms of $\mathbf{\Delta}^*$ determine the jumps of $\mathbf{X}^*_\mu$, whose law coincides with that of **X** through the following explicit construction. Given $\Delta^*$ and $n \in \mathbb{N}$, we construct $\mathbf{X}^{*n}_\mu = (X^{*n}_t, t \geq 0)$ on $2^{S \cap [n]}$ by putting

(4)
- $X^{*n}_0 = \emptyset$,
- $X^{*n}_t = X^{*n}_{t-} \triangle (\Delta_t \cap [n])$, if $(t, \Delta_t)$ is an atom of $\Delta^*$ such that $\Delta_t \cap [n] \neq \emptyset$, and
- $(X^{*n}_t, t \geq 0)$ is constant between atom times in $\Delta^*$ that affect $S \cap [n]$.

Theorem 4.6 covers the corresponding description of general combinatorial Lévy processes.

2.2. **Exchangeable processes.** Theorems 2.4 and 2.5 apply whether $S$ is finite or infinite. We obtain a more precise description of the characteristic measure $\mu$ under the additional assumption that **X** is exchangeable. In this case, the characteristic behavior of $\mu$ varies depending on whether $S$ is finite or infinite.

For processes $\mathbf{X} = (X_t, t \geq 0)$ and $\mathbf{X}' = (X'_t, t \geq 0)$, we write $\mathbf{X} =_{\mathcal{D}} \mathbf{X}'$ to denote that **X** and $\mathbf{X}'$ have the same finite-dimensional distributions, that is,

$$(X_{t_1}, \ldots, X_{t_r}) =_{\mathcal{D}} (X'_{t_1}, \ldots, X'_{t_r}) \quad \text{for all } 0 \leq t_1 \leq \cdots \leq t_r < \infty.$$

For $A \subseteq S \subseteq \mathbb{N}$ and any permutation $\sigma : S \to S$, we denote the *relabeling of $A$ by $\sigma$* by $A^\sigma$, where

$$i \in A^\sigma \quad \text{if and only if} \quad \sigma(i) \in A.$$

We call **X** *exchangeable* if $\mathbf{X} =_{\mathcal{D}} \mathbf{X}^\sigma = (X^\sigma_t, t \geq 0)$ for all permutations $\sigma : S \to S$.

By Theorem 2.4, the discrete time increments of $\mathbf{X} = (X_m, m \geq 0)$ are independent and identically distributed from a probability measure $\mu$ on $2^S$. Under the additional assumption that **X** is exchangeable, $\mu$ must also be exchangeable in the sense that $\mu(\mathcal{B}) = \mu(\mathcal{B}^\sigma)$ for all measurable $\mathcal{B} \subseteq 2^S$ and all permutations $\sigma : S \to S$, where $\mathcal{B}^\sigma = \{A^\sigma : A \in \mathcal{B}\}$.

2.2.1. *Infinite case.* The case of $S = \mathbb{N}$ is treated most directly by de Finetti's theorem [18]. Any probability measure $\nu$ on $[0, 1]$ induces an exchangeable measure $\nu^*$ on $2^{\mathbb{N}}$ by

(5) $$\nu^*(\{A^* \in 2^{\mathbb{N}} : A^* \cap [n] = A\}) = \int_{[0,1]} p^{|A|}(1-p)^{n-|A|} \nu(dp), \quad A \subseteq [n], \quad n \in \mathbb{N},$$

where $|A|$ denotes the cardinality of $A \subseteq [n]$. (Defining $\nu^*$ on sets of the form $\{A^* \in 2^{\mathbb{N}} : A^* \cap [n] = A\}$ for every $A \subseteq [n]$, $n \in \mathbb{N}$, is sufficient to uniquely determine $\nu^*$ on all of $2^{\mathbb{N}}$ because these cylinder sets are a generating $\pi$-system of the Borel $\sigma$-field on $2^{\mathbb{N}}$.) de Finetti's theorem [18] gives the converse: every exchangeable probability measure $\mu$ on $2^{\mathbb{N}}$ corresponds to a unique probability measure $\nu$ on $[0, 1]$ so that $\mu = \nu^*$, that is, $\mu$ is the $\nu$-mixture defined in (5).

In continuous time, we get the following more precise Lévy–Itô–Khintchine-type interpretation of the characteristic measure $\mu$ from Theorem 2.5.

**Theorem 2.7.** *Let $\mathbf{X} = (X_t, t \geq 0)$ be an exchangeable combinatorial Lévy process on $2^{\mathbb{N}}$. Then there exists a unique measure $\nu$ on $[0, 1]$ satisfying*

(6) $$\nu(\{0\}) = 0 \quad \text{and} \quad \int_{[0,1]} s \, \nu(ds) < \infty$$

*and a unique constant* $\mathbf{c} \geq 0$ *such that* $\mathbf{X} =_{\mathcal{D}} \mathbf{X}^*_\mu$, *the $\mu$-canonical Lévy process defined in* (4) *above with*

$$\mu = \nu^* + \mathbf{c} \sum_{i=1}^{\infty} \epsilon_i, \tag{7}$$

*where $\nu^*$ is defined as in* (5), *with $\nu$ now possibly an infinite measure, and $\epsilon_i$ is the unit mass at $\{i\} \subset \mathbb{N}$ for each $i \in \mathbb{N}$.*

We call (7) the *Lévy–Itô–Khintchine representation* for set-valued Lévy processes; see Theorem 4.18 for the corresponding theorem on general state spaces.

2.2.2. *Finite case.* When $\mathbf{X} = (X_t, t \geq 0)$ is exchangeable and evolves on $2^{[n]}$ for $n \in \mathbb{N}$, the infinitesimal jump measure $\mu$ on $2^{[n]}$ is determined instead by the finite exchangeable characterization of Diaconis and Freedman [19]. In this case, the extreme points of $2^{[n]}$ are in correspondence with the integers $0, 1, \ldots, n$. More specifically, let $\mathfrak{U}_{k:n}$ denote the law of the $\{0, 1\}$-valued sequence $(X_1, \ldots, X_n)$ obtained by recording the labels on $n$ balls sampled without replacement from an urn with $k = 0, 1, \ldots, n$ balls labeled 1 and $n - k$ balls labeled 0. Any such outcome from $\mathfrak{U}_{k:n}$ will have exactly $k$ ones and $n - k$ zeros, whence

$$\mathfrak{U}_{k:n}(x_1, \ldots, x_n) = \frac{1}{\binom{n}{k}}, \qquad \text{provided} \qquad \sum_{i=1}^{n} x_i = k.$$

It follows from the main theorem of [19] that there exists a unique $(p_0, p_1, \ldots, p_n)$ such that each $p_k \geq 0$ and $\sum_{k=0}^{n} p_k = 1$ and $\mu(\cdot) = \sum_{k=0}^{n} p_k \mathfrak{U}_{k:n}(\cdot)$.

**Theorem 2.8.** *Let $\mathbf{X} = (X_t, t \geq 0)$ be an exchangeable combinatorial Lévy process on $2^{[n]}$ for $n \in \mathbb{N}$. Then there exists a unique $(p_1, \ldots, p_n)$ satisfying $p_k \geq 0$ and $\sum_{k=1}^{n} p_k = 1$ and a unique constant $\mathbf{c} \geq 0$ such that $\mathbf{X} =_{\mathcal{D}} \mathbf{X}^*_\mu$, the $\mu$-canonical Lévy process defined above with*

$$\mu = \mathbf{c} \sum_{k=1}^{n} p_k \mathfrak{U}_{k:n}.$$

**Remark 2.9.** *We implicitly force $p_0 \equiv 0$ above in accord with the lefthand side of* (3).

2.3. **Projecting X to a space of extreme points.** Exchangeable combinatorial Lévy processes project to well behaved processes on an appropriate state space of extreme points. These extreme points index the class of measures that are ergodic with respect to the action of the symmetric group by relabeling sets.

For instance, we can project exchangeable processes $\mathbf{X} = (X_m, m \geq 0)$ on $2^\mathbb{N}$ into $[0, 1]$ by $X_m \mapsto \pi(X_m)$, where

$$\pi(X_m) := \lim_{n \to \infty} n^{-1} |X_m \cap [n]| \tag{8}$$

is the limiting frequency of elements in $X_m$, for every $m \geq 0$. By de Finetti's theorem and the strong law of large numbers, $\pi(\mathbf{X}) := (\pi(X_m), m \geq 0)$ exists almost surely whenever $\mathbf{X}$ is exchangeable. Furthermore, by independence of $X_{m-1}$ and the increments $(\Delta_r, r \geq m)$ from Theorem 2.4, we observe

$$\pi(X_m) =_{\mathcal{D}} \pi(X_{m-1})(1 - \pi(\Delta_m)) + (1 - \pi(X_{m-1}))\pi(\Delta_m),$$

so that $\pi(\mathbf{X})$ is also a Markov chain on $[0, 1]$. In continuous time, the projected process $((\pi(X_t), 1 - \pi(X_t)), t \geq 0)$ exists almost surely and exhibits the Feller property as a process on the 1-dimensional simplex equipped with the Euclidean topology; see Theorem 4.19.

In the finite case, we can also easily project any $\mathbf{X}$ on $2^{[n]}$, $n \in \mathbb{N}$, into $\{0, \ldots, n\}$ by defining

(9) $$\pi_n(X_m) = |X_m|, \quad m \geq 0,$$

the cardinality of the set $X_m$. This projection exists regardless of whether $\mathbf{X}$ is exchangeable, but in general $\pi_n(\mathbf{X}) = (\pi_n(X_m), m \geq 0)$ is a Markov chain only if $\mathbf{X}$ is exchangeable.

2.4. **Extending the set-valued case.** When moving beyond the set-valued case, the projection operation $\pi : 2^{\mathbb{N}} \to [0, 1]$ must be replaced by the more technically involved notion of a *combinatorial limit* $\|\cdot\|$, which maps a combinatorial object $M$ to an exchangeable probability measure $\|M\|$ on the space inhabited by $M$.

When $A \subseteq \mathbb{N}$, we define $\|A\|$ as follows. For any injection $\varphi : [m] \to \mathbb{N}$ and $A \subseteq \mathbb{N}$, we define $A^\varphi \subseteq [m]$ by

$$i \in A^\varphi \quad \text{if and only if} \quad \varphi(i) \in A.$$

For any $S \subseteq [m]$, we define the *limiting density of $S$ in $A$* by

$$\delta(S; A) := \lim_{n \to \infty} \frac{1}{n^{\downarrow m}} \sum_{\text{injections } \varphi:[m] \to [n]} \mathbf{1}\{A^\varphi = S\}, \quad \text{if it exists,}$$

where $n^{\downarrow m} := n(n-1) \cdots (n-m+1)$ and $\mathbf{1}\{\cdot\}$ is the indicator function of the event described by $\cdot$. (As we discuss later, existence of $\delta(S; A)$ is guaranteed whenever $A$ is the realization of an exchangeable random set.) Together the collection $(\delta(S; A), S \in \bigcup_{m \in \mathbb{N}} 2^{[m]})$ determines a unique, exchangeable probability measure $\mu$ on $2^{\mathbb{N}}$ with

$$\mu(\{A^* \in 2^{\mathbb{N}} : A^* \cap [m] = S\}) = \delta(S; A), \quad S \subseteq [m],$$

for every $m \in \mathbb{N}$. We denote this probability measure by $\|A\|$.

In the set-valued case, $\|A\|$ and $\pi(A)$ encode the same probability measure by noting that $\pi(A) = p$ implies

$$\|A\|(\{A^* \in 2^{\mathbb{N}} : A^* \cap [m] = S\}) = p^{|S|}(1-p)^{m-|S|}, \quad S \subseteq [m].$$

This equivalence is not obvious, but it follows directly from de Finetti's theorem. For general structures there is no such simplification, so we must resort to the more technical definition of $\|A\|$ in terms of the limiting densities $\delta(S; A)$, which we introduce formally in Section 3.2.

When $M$ is a combinatorial object labeled by the finite set $[n]$, the analog to the projection $\pi_n : 2^{[n]} \to \{0, 1, \ldots, n\}$ in (9) above is defined by $M \mapsto \langle M \rangle_\cong$, where $\langle M \rangle_\cong$ is the equivalence class of all objects $M'$ that are isomorphic to $M$ under relabeling, that is,

$$\langle M \rangle_\cong = \{M' : \text{there exists permutation } \sigma : [n] \to [n] \text{ such that } M^\sigma = M'\}.$$

In the set-valued case, $A \subseteq [n]$ with $|A| = m$ has

$$\langle A \rangle_\cong = \{\{i_1, \ldots, i_m\} : 1 \leq i_1 < \cdots < i_m \leq n\},$$

the collection of all $\binom{n}{m}$ size-$m$ subsets of $[n]$.

Our main theorems lift the foregoing ideas for set-valued processes to Lévy processes on arbitrary combinatorial objects. General combinatorial structures no longer have the simple 1-dimensional structure of subsets and, thus, require more care. To aid the exposition, we frame our main theorems in the context of the more tangible cases of set- and graph-valued processes whenever possible.

## 3. Combinatorial structures

The above examples are special cases of what we generally call combinatorial structures. Below we employ the usual notation $(x_1, \ldots, x_n)$ and $\{x_1, \ldots, x_n\}$ to denote ordered and unordered sets, respectively.

**Definition 3.1** (Combinatorial structures). *A signature $\mathcal{L}$ is a finite list $(i_1, \ldots, i_k)$ of nonnegative integers for which $0 \leq i_1 \leq \cdots \leq i_k$. Given a signature $\mathcal{L} = (i_1, \ldots, i_k)$ and a set $S$, a combinatorial structure with signature $\mathcal{L}$ over $S$ is a collection $M = (S; M_1, \ldots, M_k)$ such that $M_j \subseteq S^{i_j}$ for every $j = 1, \ldots, k$, with the convention $S^0 := \{\odot\}$ for $\odot$ the $S$-valued vector of length $0$. We call $i_j$ the arity of $M_j$ for each $j = 1, \ldots, k$. We alternatively call $M$ an $\mathcal{L}$-structure or simply a structure when its signature is understood. We write $\mathcal{L}_S$ to denote the set of $\mathcal{L}$-structures over $S$.*

**Remark 3.2.** *We call $M = (S; M_1, \ldots, M_k)$ a countable structure if $S$ is countable and a finite structure if $S$ is finite.*

**Remark 3.3** (Components with arity 0). *While we allow $i_1 = \cdots = i_k = 0$ in Definition 3.1, we disallow it from our main theorems. By the convention $S^0 = \{\odot\}$, the space $\mathcal{L}_S$ of structures with signature $\mathcal{L} = (0)$ consists of the two elements $(S; \emptyset)$ and $(S; \{\odot\})$. For $k > 1$, the structure $M = (S; M_1, \ldots, M_k)$ with signature $(0, \ldots, 0)$ corresponds to an element in the hypercube, which is of interest in various applications including the design of experiments but for which the labeling set $S$ plays no role. Therefore, although the case $i_1 = \cdots = i_k = 0$ is still a nontrivial state space on which to define a process, the state space is finite and, thus, lies outside the jurisdiction of our main theorems for infinite structures.*

**Remark 3.4** (Null structure). *In Section 4.3, we define the natural extension to the empty signature $\mathcal{L} = ()$. The only structure with this signature is the null structure $M = (\mathbb{N}; )$ without any relations. The significance of this structure becomes clear in Section 4.3.*

**Example 3.5** (Common examples). *In terms of Definition 3.1, a subset $A \subseteq S \subseteq \mathbb{N}$ is a combinatorial structure with $\mathcal{L} = (1)$, that is, $A \subseteq \mathbb{N}$ corresponds to $(S; A)$. A directed graph $G$ with vertex set $S$ and edge set $E \subseteq S \times S$ is a structure with $\mathcal{L} = (2)$, that is, $G = (S; E)$. (Our definition here permits self-loops in $G$.) Taking $\mathcal{L} = (1, 2)$, we obtain $M = (S; A, E)$, which corresponds to a graph $(S; E)$ and a designated subset, or community, of vertices $A \subseteq S$. For $\mathcal{L} = (1, 2, 3)$, $M = (S; M_1, M_2, M_3)$ represents first-, second-, and third-order interactions among a collection of particles or among statistical units in a designed experiment.*

**Remark 3.6** (Notation). *From now on, we use $X$ to denote random $\mathcal{L}$-structures, with $\mathbf{X}$ reserved for a family of random $\mathcal{L}$-structures, that is, a stochastic process. We use other letters, often $M$ or $A$, to represent generic (non-random) structures.*

The act of selection $S' \subseteq S$ induces a natural *restriction* operation $\mathcal{L}_S \to \mathcal{L}_{S'}$ by $M \mapsto M|_{S'}$, where

$$M|_{S'} := (S'; M_1 \cap S'^{i_1}, \ldots, M_k \cap S'^{i_k}). \tag{10}$$

Any permutation $\sigma : S \to S$ determines a *relabeling* operation $\mathcal{L}_S \to \mathcal{L}_S$ by $M \mapsto M^\sigma$, where $M^\sigma := (S; M_1^\sigma, \ldots, M_k^\sigma)$ is defined by

$$(a_1, \ldots, a_{i_j}) \in M_j^\sigma \quad \text{if and only if} \quad (\sigma(a_1), \ldots, \sigma(a_{i_j})) \in M_j \quad \text{for each } j = 1, \ldots, k. \tag{11}$$

Combining (10) and (11), we define the image of $M \in \mathcal{L}_S$ by any injection $\varphi : S' \to S$ as $M^\varphi = (S'; M_1^\varphi, \ldots, M_k^\varphi) \in \mathcal{L}_{S'}$, where

(12) $\quad (a_1, \ldots, a_{i_j}) \in M_j^\varphi \quad$ if and only if $\quad (\varphi(a_1), \ldots, \varphi(a_{i_j})) \in M_j \quad$ for each $j = 1, \ldots, k$.

Under these operations, the space $\mathcal{L}_\mathbb{N}$ of countable combinatorial structures comes furnished with the product discrete topology induced by the ultrametric

(13) $\quad d(M, M') := 1/(1 + \sup\{n \in \mathbb{N} : M|_{[n]} = M'|_{[n]}\}), \quad M, M' \in \mathcal{L}_\mathbb{N}$,

with the convention $1/\infty = 0$. Under (13), $(\mathcal{L}_\mathbb{N}, d)$ is a compact, separable, and Polish metric space, which we equip with its Borel $\sigma$-field.

When $S \subset \mathbb{N}$ is finite, the corresponding space $\mathcal{L}_S$ is also finite. In this case, we equip $\mathcal{L}_S$ with its discrete metric and we take the $\sigma$-field as the set of all subsets of $\mathcal{L}_S$.

3.1. **Combinatorial increments.** For any $S \subseteq \mathbb{N}$ and $M = (S; M_1, \ldots, M_k) \in \mathcal{L}_S$, we write

$$M_j(a) = \mathbf{1}\{a \in M_j\} := \begin{cases} 1, & a \in M_j, \\ 0, & \text{otherwise,} \end{cases}$$

for each $a = (a_1, \ldots, a_{i_j}) \in S^{i_j}$, $j = 1, \ldots, k$. We define the *increment* between $M$ and $M'$ in $\mathcal{L}_S$ by $M \triangle M' = \triangle(M, M') := (S; \Delta_1, \ldots, \Delta_k)$, where

(14) $\quad a \in \Delta_j \quad$ if and only if $\quad M_j(a) \neq M'_j(a)$,

for each $a = (a_1, \ldots, a_{i_j}) \in S^{i_j}$, $j = 1, \ldots, k$. For example, when $\mathcal{L} = (1)$, $M \triangle M'$ is the *symmetric difference* between subsets of $\mathbb{N}$ as in (1); when $\mathcal{L} = (2)$, $M \triangle M'$ is the directed graph whose edges are the pairs $(i, j)$ at which $M$ and $M'$ differ; and so on. Importantly, the increment between any two $\mathcal{L}$-structures is also an $\mathcal{L}$-structure with the same base set.

The spaces of $\mathcal{L}$-structures we consider can be regarded as a group $(\mathcal{L}_\mathbb{N}, \triangle)$, with group action given by the increment operation $\triangle$ defined above. In particular, every $M \in \mathcal{L}_S$ acts on $\mathcal{L}_S$ by $M' \mapsto M \triangle M'$. Defined in this way, $(\mathcal{L}_S, \triangle)$ is a transitive, abelian group with identity given by the empty structure $\mathbf{0}_S^\mathcal{L} := (S; \emptyset, \ldots, \emptyset)$ and for which every element $M \in \mathcal{L}_S$ is its own inverse. The group structure of $\mathcal{L}_S$ enriches $\mathcal{L}_S$ and underlies several key properties of combinatorial Lévy processes. Furthermore, $\mathcal{L}_S$ is partially ordered and has minimum element $\mathbf{0}_S^\mathcal{L}$ under pointwise inclusion, that is, $M \leq M'$ if and only if $M_j(a) \leq M'_j(a)$ for every $a \in S^{i_j}$, for all $j = 1, \ldots, k$.

3.2. **Exchangeability and combinatorial limits.** de Finetti's theorem, Diaconis and Freedman's theorem, and the Aldous–Hoover theorem permit the study of infinite and finite exchangeable sequences and graphs by projecting into a limit space, for example, the unit interval, an initial segment of the non-negative integers, or the space of graph limits. We observe analogous behavior in general.

3.2.1. *Countable structures.* The example in Section 2 shows that much of the structural behavior of an exchangeable set-valued Lévy process is determined by its projection into the unit interval. Our main theorems extend this idea to describe exchangeable combinatorial Lévy processes through their induced behavior in the appropriate limit space.

Infinite exchangeable combinatorial $\mathcal{L}$-structures admit a representation in a space of *combinatorial limits*. As mentioned in Section 2, the combinatorial limit of $M \in \mathcal{L}_\mathbb{N}$ cannot be described as simply as the projection of $A \subseteq \mathbb{N}$ to its limiting frequency $\pi(A)$ as in (8). To see why, consider $A, A' \subseteq \mathbb{N}$ and let $M = (\mathbb{N}; A, A')$ be the associated $(1, 1)$-structure. Although $M$ consists of a pair of subsets, the individual frequencies $\pi(A)$ and $\pi(A')$ are not sufficient to summarize the full structure of $M$: if we construct $A$ by including each element $i \in \mathbb{N}$ independently with probability $p \in (0, 1)$ and then defining $A' = A$, then $\pi(A) = \pi(A') = p$ with probability 1; but if we define $A$ and $A'$ as independent and identically distributed so that each element has probability $p \in (0, 1)$ of appearing in $A$, respectively $A'$, then $\pi(A) = \pi(A') = p$ with probability 1, but $\mathbb{P}\{A = A'\} = 0$. In both cases, $(\pi(A), \pi(A')) = (p, p)$ with probability 1, but the structure of $M = (\mathbb{N}; A, A')$ is vastly different under the two constructions. The pair $(\pi(A), \pi(A'))$ does not capture all relevant structural features of $(\mathbb{N}; A, A')$, motivating the following definition.

**Definition 3.7** (Homomorphism density). *For any signature $\mathcal{L}$ and finite subsets $S' \subseteq S \subset \mathbb{N}$, we define the* homomorphism density *of $A \in \mathcal{L}_{S'}$ in $M \in \mathcal{L}_S$ by*

(15) $$\delta(A; M) := \frac{1}{|S|^{\downarrow |S'|}} \sum_{\varphi: S' \to S} \mathbf{1}\{M^\varphi = A\},$$

*where the sum is over injections $\varphi : S' \to S$, $|S|$ denotes the cardinality of $S \subseteq \mathbb{N}$, and $n^{\downarrow m} := n(n-1) \cdots (n-m+1)$ is the falling factorial function. For brevity, we refer to (15) as the* density *of $A$ in $M$.*

Intuitively, $\delta(A; M)$ is the probability that $M^\varphi = A$ for $\varphi$ chosen uniformly at random among all injections $S' \to S$. For fixed $M \in \mathcal{L}_S$, the density function $\delta(\cdot; M)$ determines a probability measure on $\mathcal{L}_{S'}$ for every $S' \subseteq S$. For $M \in \mathcal{L}_\mathbb{N}$ and $A \in \mathcal{L}_{[m]}$, we define the *limiting density* of $A$ in $M$ by

(16) $$\delta(A; M) := \lim_{n \to \infty} \delta(A; M|_{[n]}), \quad \text{if it exists.}$$

Provided each of the limits $\delta(A; M)$, $A \in \mathcal{L}_{[n]}$, exists, the collection of homomorphism densities $(\delta(A; M), A \in \mathcal{L}_{[n]})$ determines a probability measure on $\mathcal{L}_{[n]}$ by the bounded convergence theorem. If (16) exists for every $A \in \bigcup_{n \in \mathbb{N}} \mathcal{L}_{[n]}$, then the family of distributions $((\delta(A; M), A \in \mathcal{L}_{[n]}), n \in \mathbb{N})$ determines a unique probability measure on $\mathcal{L}_\mathbb{N}$, which we denote by $\|M\|$.

**Definition 3.8** (Combinatorial limit). *The* combinatorial limit *$\|M\|$ of $M \in \mathcal{L}_\mathbb{N}$ is the unique probability measure $\mu$ on $\mathcal{L}_\mathbb{N}$ such that*

(17) $$\mu(\{A^* \in \mathcal{L}_\mathbb{N} : A^*|_{[m]} = A\}) = \delta(A; M), \quad A \in \mathcal{L}_{[m]}, \quad m \in \mathbb{N},$$

*provided the limit $\delta(A; M)$ exists for every $A \in \bigcup_{m \in \mathbb{N}} \mathcal{L}_{[m]}$. For brevity we often write*

$$\|M\|(A) := \|M\|(\{A^* \in \mathcal{L}_\mathbb{N} : A^*|_{[m]} = A\}) \quad \text{for each } A \in \mathcal{L}_{[m]}, \quad m \in \mathbb{N}.$$

Lovász and Szegedy [34] defined the concept of a *graph limit* in terms of the limiting homomorphism densities of all finite subgraphs within a sequence of graphs. Definition 3.8 extends the Lovász–Szegedy notion to the more general setting of combinatorial structures from Definition 3.1. The space of exchangeable, dissociated probability measures plays an important role in the study of exchangeable structures.

**Definition 3.9** (Exchangeable and dissociated $\mathcal{L}$-structures). *For any $S \subseteq \mathbb{N}$ and any signature $\mathcal{L}$, a random $\mathcal{L}$-structure $X$ over $S$ is*
- exchangeable, *if $X^\sigma =_\mathcal{D} X$ for all permutations $\sigma : S \to S$, and*
- dissociated, *if $X|_T$ and $X|_{T'}$ are independent whenever $T, T' \subseteq S$ are disjoint.*

12 HARRY CRANEWhen $A \subseteq \mathbb{N}$ is a random subset, exchangeable and dissociated corresponds to each $i \in \mathbb{N}$ being present in $A$ independently with the same probability, which explains why the projection $\pi(A)$ in (8) is enough to determine the combinatorial limit $\|X\|$ of the random (1)-structure $X = (\mathbb{N}; A)$; see Equation (8) and the discussion at the end of Section 2. In Proposition 6.2, we prove that the combinatorial limit of any exchangeable $\mathcal{L}$-structure exists with probability 1.

**Definition 3.10** (Combinatorial limit space). *For any signature $\mathcal{L}$, we write $\mathcal{E}_\mathcal{L}$ to denote the space of exchangeable, dissociated probability measures on $\mathcal{L}_\mathbb{N}$.*

**Remark 3.11.** *The notation $\mathcal{E}_\mathcal{L}$ should be understood to stand for the space of* ergodic *measures with respect to the action of the symmetric group.*

As every $W \in \mathcal{E}_\mathcal{L}$ is a probability measure on $\mathcal{L}_\mathbb{N}$, we write $W(A)$, $A \in \mathcal{L}_{[n]}$, as shorthand for
$$W(A) := W(\{A^* \in \mathcal{L}_\mathbb{N} : A^*|_{[n]} = A\}), \quad A \in \mathcal{L}_{[n]}.$$
We then define the distance between $W, W' \in \mathcal{E}_\mathcal{L}$ by
$$(18) \qquad d(W, W') = \sum_{n \in \mathbb{N}} 2^{-n} \sum_{A \in \mathcal{L}_{[n]}} |W(A) - W'(A)|.$$

Under (18), $\mathcal{E}_\mathcal{L}$ is a closed subset of the space of all probability measures on $\mathcal{L}_\mathbb{N}$ and, therefore, is compact. We define $\|M\| = \partial$ whenever at least one of the limiting densities $\delta(A; M)$ does not exist. With this, we define $d(W, \partial) = 2$ for all $W \in \mathcal{E}_\mathcal{L}$ and equip $\mathcal{E}_\mathcal{L}$ with the Borel $\sigma$-field induced by this metric.

3.2.2. *Finite structures.* For $S \subset \mathbb{N}$ finite, $\mathcal{L}_S$ is also finite and we need not take limits to encode the information contained in (15). Rather, we define $\mathcal{UL}_S$ as the quotient space of $\mathcal{L}_S$ under the equivalence relation
$$(19) \qquad M \cong M' \quad \text{if and only if} \quad \text{there exists } \sigma : S \to S \text{ such that } M^\sigma = M',$$
writing specifically $\langle M \rangle_\cong \in \mathcal{UL}_S$ to denote the equivalence class of $M$ under $\cong$. For $\mathcal{L} = (1)$ and $S = [n]$, $\mathcal{UL}_{[n]}$ partitions $\mathcal{L}_{[n]}$ into subsets of cardinality $0, 1, \ldots, n$, in agreement with the discussion surrounding Theorem 2.8 and (9).

For any signature $\mathcal{L}$, $n \in \mathbb{N}$, and $Y \in \mathcal{UL}_{[n]}$, we define $\mathfrak{U}_Y$ to be the uniform distribution on $\{M \in \mathcal{L}_{[n]} : \langle M \rangle_\cong = Y\}$.

**Proposition 3.12.** *Let $\mathcal{L}$ be any signature, $n \in \mathbb{N}$, and $\mu$ be any exchangeable probability measure on $\mathcal{L}_{[n]}$. Then there exists a unique $(p_Y)_{Y \in \mathcal{UL}_{[n]}}$ such that $p_Y \geq 0$ for all $Y \in \mathcal{UL}_{[n]}$, $\sum_{Y \in \mathcal{UL}_{[n]}} p_Y = 1$, and*
$$(20) \qquad \mu(\cdot) = \sum_{Y \in \mathcal{UL}_{[n]}} p_Y \mathfrak{U}_Y(\cdot).$$

*Proof.* Suppose $X \sim \mu$. Then exchangeability of $\mu$ implies
$$\mathbb{P}\{X = M \mid \langle X \rangle_\cong\} = \mathbb{P}\{X = M^\sigma \mid \langle X \rangle_\cong\} \quad \text{for all permutations } \sigma : [n] \to [n].$$
By (19), $\langle M \rangle_\cong = \langle M^\sigma \rangle_\cong$ for all permutations $\sigma : [n] \to [n]$. Now, define $p_Y := \mathbb{P}\{\langle X \rangle_\cong = Y\}$ and note that
$$\mathfrak{U}_Y(M) = \begin{cases} |\{M' \in \mathcal{L}_{[n]} : \langle M' \rangle_\cong = Y\}|^{-1}, & \langle M \rangle_\cong = Y, \\ 0, & \text{otherwise.} \end{cases}$$

By the law of total probability,

$$\begin{aligned} \mu(M) &= \sum_{Y \in \mathcal{UL}_{[n]}} \mathbb{P}\{X = M \mid \langle X \rangle_{\cong} = Y\} \mathbb{P}\{\langle X \rangle_{\cong} = Y\} \\ &= \sum_{Y \in \mathcal{UL}_{[n]}} \mathfrak{U}_Y(M) \mathbb{P}\{\langle X \rangle_{\cong} = Y\} \\ &= \sum_{Y \in \mathcal{UL}_{[n]}} p_Y \mathfrak{U}_Y(M), \quad \text{for all } M \in \mathcal{L}_{[n]}. \end{aligned}$$

□

## 4. Summary of main theorems

**Remark 4.1.** *All theorems below assume a signature $\mathcal{L} = (i_1, \ldots, i_k)$ for which $i_k \geq 1$.*

### 4.1. General combinatorial Lévy processes.
Recall the definition of the increment $\triangle : \mathcal{L}_S \times \mathcal{L}_S \to \mathcal{L}_S$ in (14) and $\mathbf{0}_S^{\mathcal{L}} = (S; \emptyset, \ldots, \emptyset)$.

**Definition 4.2** (Combinatorial Lévy process). *For any signature $\mathcal{L}$ and $S \subseteq \mathbb{N}$, we call $\mathbf{X} = (X_t, t \geq 0)$ on $\mathcal{L}_S$ a combinatorial Lévy process if it has*

- $X_0 = \mathbf{0}_S^{\mathcal{L}}$,
- *stationary increments, that is, $\triangle(X_{t+s}, X_s) =_{\mathcal{D}} X_t$ for all $s, t \geq 0$,*
- *independent increments, that is, $\triangle(X_{t_1}, X_{t_0}), \ldots, \triangle(X_{t_k}, X_{t_{k-1}})$ are independent for all $0 \leq t_0 \leq t_1 \leq \cdots \leq t_k < \infty$, and*
- *càdlàg sample paths, that is, $t \mapsto X_t$ is right continuous and has left limits under the product discrete topology induced by (13).*

**Remark 4.3.** *The first condition above, $X_0 = \mathbf{0}_S^{\mathcal{L}}$, is akin to the condition $X_0 = 0$ for $\mathbb{R}$-valued Lévy processes. By stationarity and independence of increments, there is no loss of generality in assuming $X_0 = \mathbf{0}_S^{\mathcal{L}}$. A combinatorial Lévy process $\mathbf{X}^x$ with initial state $X_0 = x$ can be obtained from $\mathbf{X} = (X_t, t \geq 0)$ started at $\mathbf{0}_S^{\mathcal{L}}$ by putting $\mathbf{X}^x := (X_t^x, t \geq 0)$ with $X_t^x = X_t \triangle x$ for all $t \geq 0$.*

In discrete time, combinatorial Lévy processes are analogous to random walks. Most of their structural properties follow directly from Definition 4.2.

**Definition 4.4** (Combinatorial random walk). *A (combinatorial) random walk on $\mathcal{L}_S$ with increment distribution $\mu$ and initial state $x$ is a discrete time process $\mathbf{X}^x = (X_m, m \geq 0)$ with $X_0 = x$ and*

$$(21) \qquad X_m =_{\mathcal{D}} \triangle(X_{m-1}, \Delta_m), \quad m \geq 1,$$

*where $\Delta_1, \Delta_2, \ldots$ are i.i.d. from $\mu$.*

**Theorem 4.5.** *Let $\mathbf{X} = (X_m, m \geq 0)$ be a discrete time combinatorial Lévy process on $\mathcal{L}_S$. Then there exists a unique probability measure $\mu$ on $\mathcal{L}_S$ such that $\mathbf{X} =_{\mathcal{D}} \mathbf{X}_\mu^* = (X_m^*, m \geq 0)$, where $\mathbf{X}_\mu^*$ is a combinatorial random walk on $\mathcal{L}_S$ with initial state $X_0 = \mathbf{0}_S^{\mathcal{L}}$ and increment distribution $\mu$.*

In continuous time, a combinatorial Lévy process on $\mathcal{L}_\mathbb{N}$ must balance its behavior so that its sample paths satisfy the càdlàg requirement. Since each $\mathcal{L}_{[n]}$ is a finite state space, $\mathbf{X}^{[n]} := (X_t|_{[n]}, t \geq 0)$ can jump only finitely often in bounded time intervals. On the other hand, since we have ruled out the case $i_1 = \cdots = i_k = 0$, $\mathbf{X} = (X_t, t \geq 0)$ evolves on an uncountable state space and is defined at an uncountable set of times; therefore, $\mathbf{X}$ can experience possibly infinitely many discontinuities in bounded time intervals. Condition (22) in Theorem 4.6 strikes the balance.

**Theorem 4.6.** *Let* $\mathbf{X} = (X_t, t \geq 0)$ *be a continuous time combinatorial Lévy process on* $\mathcal{L}_S$. *Then there is a unique measure $\mu$ on $\mathcal{L}_S$ satisfying*

$$\text{(22)} \qquad \mu(\{\mathbf{0}_S^{\mathcal{L}}\}) = 0 \quad \text{and} \quad \mu(\{M^* \in \mathcal{L}_S : M^*|_{S \cap [n]} \neq \mathbf{0}_{S \cap [n]}^{\mathcal{L}}\}) < \infty \text{ for all } n \in \mathbb{N}$$

*such that the infinitesimal jump rates of $\mathbf{X}$ satisfy*

$$\text{(23)} \qquad \lim_{t \downarrow 0} \frac{1}{t} \mathbb{P}\{X_t \in d\Delta\} = \mu(d\Delta), \quad \Delta \in \mathcal{L}_S \setminus \{\mathbf{0}_S^{\mathcal{L}}\},$$

*where convergence in* (23) *is understood in the sense of vague convergence of $\sigma$-finite measures.*

**Remark 4.7.** *Condition* (22) *applies whether $S$ is finite or infinite. If $S \subseteq \mathbb{N}$ is finite, then $S \cap [n] = S$ for all large $n \in \mathbb{N}$ and the righthand side of* (22) *simply means that $\mu$ is a finite measure.*

The limit in (23) is well defined on account of the Feller property for combinatorial Lévy processes. The stationary and independent increments assumptions imply that $\mathbf{X}$ is a time homogeneous Markov process with transition law determined by the Markov semigroup $\mathbf{Q} = (\mathbf{Q}_t, t \geq 0)$, where

$$\text{(24)} \qquad \mathbf{Q}_t g(M) := \mathbb{E}g(X_t \triangle M), \quad t \geq 0,$$

for all bounded, continuous functions $g : \mathcal{L}_S \to \mathbb{R}$ and all $M \in \mathcal{L}_S$.

**Definition 4.8** (Feller property). *Let $\mathbf{Z} = (Z_t, t \geq 0)$ be a Markov process on a topological space $\mathcal{Z}$ with semigroup $\mathbf{T} = (\mathbf{T}_t, t \geq 0)$. We call $\mathbf{T}$ a* Feller semigroup *and say that $\mathbf{Z}$ has the* Feller property *if*
- $\lim_{t \downarrow 0} \mathbf{T}_t g(z) = g(z)$ *for all $z \in \mathcal{Z}$ and*
- $z \mapsto \mathbf{T}_t g(z)$ *is continuous for every $t > 0$,*

*for all bounded, continuous $g : \mathcal{Z} \to \mathbb{R}$.*

**Proposition 4.9.** *An $\mathcal{L}_S$-valued process $\mathbf{X} = (X_t, t \geq 0)$ is a combinatorial Lévy process if and only if $\mathbf{X}^{S'} = (X_t|_{S'}, t \geq 0)$ is a combinatorial Lévy process on $\mathcal{L}_{S'}$ for every $S' \subseteq S$. Moreover, an $\mathcal{L}_\mathbb{N}$-valued process $\mathbf{X}$ is a combinatorial Lévy process if and only if $\mathbf{X}^{[n]} = (X_t|_{[n]}, t \geq 0)$ is a combinatorial Lévy process on $\mathcal{L}_{[n]}$ for every $n = 1, 2, \ldots$.*

We immediately deduce the Feller property for combinatorial Lévy processes.

**Corollary 4.10.** *Every combinatorial Lévy process has the Feller property.*

**Definition 4.11** ($\sigma$-finite measures). *A measure $\mu$ on $\mathcal{L}_S$ is $\sigma$-finite if it satisfies* (22).

Given a $\sigma$-finite measure $\mu$ on $\mathcal{L}_S$, we construct the $\mu$-*canonical Lévy process* $\mathbf{X}_\mu^* = (X_t^*, t \geq 0)$ from a Poisson point process $\Delta^* = \{(t, \Delta_t^*)\} \subseteq [0, \infty) \times \mathcal{L}_S$ with intensity measure $dt \otimes \mu$. For each $n \in \mathbb{N}$, we construct $\mathbf{X}_\mu^{*n} = (X_t^{*n}, t \geq 0)$ on $\mathcal{L}_{S \cap [n]}$ by putting

- $X_0^{*n} = \mathbf{0}_{S \cap [n]}^{\mathcal{L}}$,
- $X_t^{*n} = X_{t-}^{*n} \triangle \Delta_t^*|_{S \cap [n]}$, if $(t, \Delta_t^*) \in \Delta^*$ and $\Delta_t^*|_{S \cap [n]} \neq \mathbf{0}_{S \cap [n]}^{\mathcal{L}}$,
- and otherwise $\mathbf{X}_\mu^{*n}$ is constructed to be constant between atom times of $\Delta^*$ that affect $S \cap [n]$.

Since we construct each $\mathbf{X}_\mu^{*n}$ from the same Poisson point process $\mathbf{\Delta}^*$, the collection $(\mathbf{X}_\mu^{*n}, n \in \mathbb{N})$ is mutually compatible, that is, $\mathbf{X}_\mu^{*n}|_{S \cap [m]} := (X_t^{*n}|_{S \cap [m]}, t \geq 0) = \mathbf{X}_\mu^{*m}$ for every $m \leq n$, and, thus, determines a unique process $\mathbf{X}_\mu^* = (X_t^*, t \geq 0)$ on $\mathcal{L}_S$.

**Remark 4.12.** *Note that the above construction is necessary only when $S$ is infinite. When $S$ is finite, (22) implies that $\mu$ is a finite measure and the canonical Lévy process $\mathbf{X}_\mu^*$ can be constructed on $\mathcal{L}_S$ directly by a single pass through the above construction at any level $n \geq \max S$.*

**Theorem 4.13.** *Let $\mathbf{X}$ be a combinatorial Lévy process on $\mathcal{L}_S$ with rate measure $\mu$ as in (22). Then $\mathbf{X} =_\mathcal{D} \mathbf{X}_\mu^*$, where $\mathbf{X}_\mu^*$ is a $\mu$-canonical Lévy process. Moreover, every combinatorial Lévy process $\mathbf{X}$ has the same finite-dimensional distributions as some $\mu$-canonical Lévy process corresponding to a $\sigma$-finite measure $\mu$ on $\mathcal{L}_S$.*

4.2. **Exchangeable processes.** For any permutation $\sigma: S \to S$, we write $\mathbf{X}^\sigma := (X_t^\sigma, t \geq 0)$ to denote the image of $\mathbf{X}$ under relabeling by $\sigma$.

**Definition 4.14** (Exchangeable Lévy process). *An $\mathcal{L}_S$-valued process $\mathbf{X} = (X_t, t \geq 0)$ is exchangeable if $\mathbf{X} =_\mathcal{D} \mathbf{X}^\sigma$ for all permutations $\sigma: S \to S$.*

Definition 3.8 extends the notion of limits of large graphs from [34] to general $\mathcal{L}$-structures. For every signature $\mathcal{L}$, the limit space $\mathcal{E}_\mathcal{L}$ consists of exchangeable, dissociated probability measures on $\mathcal{L}_\mathbb{N}$. Given a measure $\nu$ on $\mathcal{E}_\mathcal{L}$, we write $\nu^*$ to denote the exchangeable measure it induces on $\mathcal{L}_\mathbb{N}$ by

$$\nu^*(\mathcal{B}) := \int_{\mathcal{E}_\mathcal{L}} W(\mathcal{B}) \nu(dW), \quad \mathcal{B} \subseteq \mathcal{L}_\mathbb{N}, \tag{25}$$

where $W(\mathcal{B})$ is the measure assigned to $\mathcal{B} \subseteq \mathcal{L}_\mathbb{N}$ by $W$. As long as $\nu$ is a probability measure on $\mathcal{E}_\mathcal{L}$, $\nu^*$ is a probability measure on $\mathcal{L}_\mathbb{N}$, but the definition in (25) is well defined for arbitrary positive measures $\nu$. When $\nu$ is a probability measure, (25) has the interpretation of first drawing $W \sim \nu$ and, given $W$, sampling a random $\mathcal{L}$-structure from $W$.

For $n \in \mathbb{N}$ and $\mathbf{p} = (p_Y)_{Y \in \mathcal{UL}_{[n]}}$, we write $\mathbf{p}^*$ to denote the exchangeable measure induced on $\mathcal{L}_{[n]}$ by

$$\mathbf{p}^*(M) = \sum_{Y \in \mathcal{UL}_{[n]}} p_Y \mathfrak{U}_Y(M), \quad M \in \mathcal{L}_{[n]}. \tag{26}$$

For any combinatorial Lévy process $\mathbf{X} = (X_t, t \geq 0)$ on $\mathcal{L}_\mathbb{N}$, we write $\|\mathbf{X}\| = (\|X_t\|, t \geq 0)$ to denote its projection into $\mathcal{E}_\mathcal{L}$, if it exists. And for any combinatorial Lévy process $\mathbf{X} = (X_t, t \geq 0)$ on $\mathcal{L}_{[n]}$, $n \in \mathbb{N}$, we write $\langle \mathbf{X} \rangle_\cong = (\langle X_t \rangle_\cong, t \geq 0)$ to denote its projection into $\mathcal{UL}_{[n]}$, which always exists. The next theorem says that $\|\mathbf{X}\|$ exists with probability 1 for infinite exchangeable combinatorial Lévy processes in discrete time. Theorem 4.19 gives the corresponding statement for continuous time combinatorial Lévy processes.

**Theorem 4.15.** (a) **Infinite case**: *Let $\mathbf{X} = (X_m, m \geq 0)$ be a discrete time exchangeable combinatorial Lévy process on $\mathcal{L}_\mathbb{N}$. Then there exists a unique probability measure $\nu$ on $\mathcal{E}_\mathcal{L}$ such that the increments of $\mathbf{X}$ are independent and identically distributed according to $\nu^*$. Moreover, the projection $\|\mathbf{X}\| = (\|X_m\|, m \geq 0)$ exists almost surely and is a Markov chain on $\mathcal{E}_\mathcal{L}$.*
   (b) **Finite case**: *Let $\mathbf{X} = (X_m, m \geq 0)$ be a discrete time exchangeable combinatorial Lévy process on $\mathcal{L}_{[n]}$, for some $n \in \mathbb{N}$. Then there exists a unique $\mathbf{p} = (p_Y)_{Y \in \mathcal{UL}_{[n]}}$ such that each $p_Y \geq 0$, $\sum_{Y \in \mathcal{UL}_{[n]}} p_Y = 1$, and the increments of $\mathbf{X}$ are independent and identically* *distributed according to* $\mathbf{p}^*$. *Moreover, the projection* $\langle \mathbf{X} \rangle_\cong = (\langle X_m \rangle_\cong, m \geq 0)$ *is a Markov chain on* $\mathcal{UL}_{[n]}$.

4.3. **Lévy–Itô structure.** The final theorems apply exclusively to countable structures. Here we often deal with unordered multisets, which we write as $s = \{s_1^{m_1}, \ldots, s_r^{m_r}\}$ with $s_1 < \cdots < s_r$ and each element $s_i$ appearing with *multiplicity* $m_i$ in $s$. With this notation, a subset $s = \{s_1, \ldots, s_r\} \subseteq \mathbb{N}$ corresponds to $\{s_1^1, \ldots, s_r^1\}$ with multiplicities omitted, allowing us to extend the notation $s \subseteq \mathbb{N}$ to indicate that $s$ is a multiset of $\mathbb{N}$. We write $|s| = \sum_{1 \leq i \leq r} m_i$ to denote the cardinality of $s = \{s_1^{m_1}, \ldots, s_r^{m_r}\}$ counted with multiplicity and $\{s\} := \{s_1, \ldots, s_r\}$ to denote the set of distinct elements in $s$. For example, $s = \{1^2, 4^3\}$ corresponds to the multiset $\{1, 1, 4, 4, 4\}$, for which $\{s\} = \{1, 4\}$. Given two multisets $s, s'$, we write $s \preceq s'$ to denote that $\{s\} \subseteq \{s'\}$ and each element appears in $s$ with multiplicity no greater than its multiplicity in $s'$. We define the intersection $s \cap s'$ to be the multiset with each element of $\{s\} \cap \{s'\}$ appearing with multiplicity equal to its minimum multiplicity in $s$ and $s'$. For example, $\{1^2, 4^3\} \preceq \{1^3, 3^1, 4^3\}$ but $\{1^2, 4^3\} \not\preceq \{1^3, 3^1, 4^2\}$, and for $s = \{1^2, 4^3\}$ and $s' = \{1^3, 2^1, 3^4, 4^2\}$ we have $s \cap s' = \{1^2, 4^2\}$. We apply the same notation for ordered multisets $a = (a_1, \ldots, a_n)$ when the order does not matter, that is, we write $s \preceq a$ to denote that the relation holds with $a$ regarded as the multiset determined by its components. For example, $a = (1, 3, 1, 4, 1, 4)$ determines the multiset $\{1^3, 3^1, 4^2\}$, for which we have $\{1^2, 4^2\} \preceq a$.

Let $\mathcal{L} = (i_1, \ldots, i_k)$ be a signature and $s = \{s_1^{m_1}, \ldots, s_r^{m_r}\} \subset \mathbb{N}$ be a multiset with $|s| = q$ for some $q = 0, 1, \ldots, i_k$. For any $M \in \mathcal{L}_\mathbb{N}$, we define the *s-substructure* $M_s^* = (\mathbb{N}; M_{s,1}^*, \ldots, M_{s,k}^*)$ as the $\mathcal{L}$-structure with

$$
(27) \qquad M_{s,j}^*(a) = \begin{cases} M_j(a), & |a| \leq |s|, \ a \preceq s, \ \{a\} = \{s\}, \\ M_j(a), & |a| > |s|, \ s \preceq a, \\ 0, & \text{otherwise,} \end{cases} \qquad a \in \mathbb{N}^{i_j}, \quad j = 1, \ldots, k.
$$

Therefore, $M_s^*$ is the $\mathcal{L}$-structure that corresponds to $M$ on supersets of $s$ and to $\mathbf{0}_\mathbb{N}^\mathcal{L}$ otherwise. As a special case, we point out that $M_\emptyset^* = M$ for all $M \in \mathcal{L}_\mathbb{N}$.

The first two separate conditions in (27) are needed to fully capture the possible behaviors in our main theorem below. In (27) the multiset $s = \{s_1^{m_1}, \ldots, s_r^{m_r}\}$ represents the elements indexing the chosen substructure of $M$. If $|s| \geq i_j$ for some component $j = 1, \ldots, k$ of $\mathcal{L} = (i_1, \ldots, i_k)$, then $M_{s,j}^*(a)$ can be nonzero only if $\{a\}$ and $\{s\}$ coincide as sets without multiplicity. If $|s| < i_j$, then $M_{s,j}^*(a)$ can be nonzero only if all elements of $s$ appear in $a$ with multiplicity at least their multiplicity in $s$. Some examples clarify this definition.

**Example 4.16.** *Let* $\mathcal{L} = (1, 2)$ *so that* $M = (\mathbb{N}; A, E)$ *is a set* $A \subseteq \mathbb{N}$ *together with a graph* $(\mathbb{N}; E)$. *For* $s = \{1^1\}$, $M_s^*$ *retains only relations in $M$ involving element 1. Specifically,* $M_{s,j}^*(a) = 0$ *for all tuples* $a = (a_{i_1}, \ldots, a_{i_j})$ *except possibly those containing element 1:*

$$M_{s,1}^*((i)) = \begin{cases} M_1((1)), & i = 1, \\ 0, & \text{otherwise,} \end{cases} \quad \text{and}$$

$$M_{s,2}^*((i, i')) = \begin{cases} M_2((i, i')), & i = 1 \text{ or } i' = 1, \\ 0, & \text{otherwise.} \end{cases}$$

*With regard to (27), $s = \{1^1\}$ implies $|s| = 1$, so that $M_{s,1}^*(a)$ is determined by the top line since $|a| = 1 = |s|$, while $M_{s,2}^*(a)$ is determined by the second line since $|a| = 2 > 1 = |s|$. Since the top line*

*requires $\{a\} = \{s\} = 1$, the only nontrivial contribution to $M^*_{s,1}$ comes from $a = (1)$. The second line requires only $s \leq a$, allowing for all $a = (i, i')$ such that $\{1\} \subseteq \{i, i'\}$.*

*We note the difference when $s = \{1^2\}$, which also has $\{s\} = \{1\}$ but should not be confused with $\{1^1\}$ in the context of (27). For $s = \{1^2\}$,*

$$M^*_{s,1}((i)) = \begin{cases} M_1((1)), & i = 1, \\ 0, & \text{otherwise,} \end{cases} \quad \text{and}$$

$$M^*_{s,2}((i, i')) = \begin{cases} M_2((1, 1)), & (i, i') = (1, 1), \\ 0, & \text{otherwise.} \end{cases}$$

*Once again, $M^*_{s,1}(a)$ is determined by the top line of (27), for which the only nontrivial contribution must have $\{a\} = \{1\}$ and, therefore, $a = (1)$. In contrast to the case $s = \{1^1\}$, however, the top line of (27) also applies to $M^*_{s,2}(a)$, since we now have $|a| = 2 = |s|$. The contribution $M^*_{s,2}(a)$ is nontrivial only if $\{a\} = \{s\} = \{1\}$, that is, $a = (1, 1)$. Therefore, $M^*_{\{1^1\}}$ and $M^*_{\{1^2\}}$ are different structures in general.*

*Finally, consider $s = \{1^1, 2^1\}$, then*

$$M^*_{s,1}((i)) = 0 \quad \text{for all } i \in \mathbb{N}$$

*since it is impossible for $1 = |a| < |s| = 2$ and $\{a\} = \{1, 2\}$. On the other hand,*

$$M^*_{s,2}((i, i')) = \begin{cases} M_2((i, i')), & (i, i') = (1, 2) \text{ or } (i, i') = (2, 1), \\ 0, & \text{otherwise.} \end{cases}$$

Any multiset $s = \{s_1^{m_1}, \ldots, s_r^{m_r}\} \subset \mathbb{N}$ determines a partition of the integer $q = |s|$, written as $\lambda(s) = (m_i^\downarrow, 1 \leq i \leq r)$, the multiplicities of $s$ listed in nonincreasing order. In general we write $\lambda \vdash q$ to indicate that $\lambda = (\lambda_1, \ldots, \lambda_r)$ is a partition of the integer $q$, which must satisfy $\lambda_1 \geq \cdots \geq \lambda_r \geq 1$ and $\lambda_1 + \cdots + \lambda_r = q$. We often omit parentheses and write $\lambda = (\lambda_1, \ldots, \lambda_r) = \lambda_1 \lambda_2 \cdots \lambda_r$.

For $j = 1, \ldots, k$ and $s = \{s_1^{m_1}, \ldots, s_r^{m_r}\}$ with $|s| = q$ for some $q = 0, 1, \ldots, i_k$, we can express each component $M^*_{s,j}$ of $M^*_s$ as a structure in and of itself with signature $(i_j - q)^{k_j} = (i_j - q, \ldots, i_j - q)$ with $k_j$ equal arities, where

$$k_j := \binom{i_j}{i_j - q} \frac{q!}{\prod_{l=1}^r \lambda_l!}$$

for $\lambda(s) = (\lambda_1, \ldots, \lambda_r) \vdash q$. If $i_j < q$, then $k_j = 0$ and we define $M^*_{s,j}$ as the null structure $(\mathbb{N}; )$ with empty signature $\mathcal{L} = ()$, as in Remark 3.4. If $i_j = q$, then $M^*_{s,j}$ has the signature $0^{k_j}$.

Note that $k_j$ is the number of all possible ways to insert the elements of $s$ in an $i_j$-tuple in any possible order. For example, $i_j = 3$ and $s = \{1^2\}$ give $q = 2$, $\lambda(s) = 2$, and $k_j = 3$, which corresponds to the three tuples

$$(*, 1, 1), \quad (1, *, 1), \quad (1, 1, *),$$

where entries $*$ can be filled with arbitrary indices. On the other hand, $i_j = 3$ and $s = \{1^1, 2^1\}$ give $q = 2$, $\lambda(s) = (1, 1)$, and $k_j = 6$ corresponding to the six tuples of the form

$$(*, 1, 2), \ (*, 2, 1), \ (1, *, 2), \ (2, *, 1), \ (1, 2, *), \ (2, 1, *).$$

In this case, we can encode $M^*_{s,j}$ by another structure $(\mathbb{N}; M^*_{s,j,1}, \ldots, M^*_{s,j,6})$, where each $M^*_{s,j,l} \subseteq \mathbb{N}^{i_j-q} = \mathbb{N}$. With the indices $l = 1, \ldots, 6$ corresponding to the ordering of tuples

above, we have

$$M^*_{s,j,1}((a)) = M_j((a, 1, 2)),$$
$$M^*_{s,j,2}((a)) = M_j((a, 2, 1)),$$
$$M^*_{s,j,3}((a)) = M_j((1, a, 2)),$$

and so on, for each $a \in \mathbb{N}^{i_j-q} = \mathbb{N}$. For $s \subset \mathbb{N}$ with $\lambda(s) = \lambda$, we write $\mathcal{L}^\lambda$ to denote the signature $((i_j - q)^{k_j} : j = 1, \ldots, k)$ of $M^*_s$ when each $M^*_{s,j}$ is encoded as an $(i_j - q)^{k_j}$-structure with any chosen convention of the ordering of tuples.

For any $\lambda = (\lambda_1, \ldots, \lambda_r) \vdash q$, we define the *canonical $\lambda$-multiset* by $s_\lambda = \{1^{\lambda_1}, 2^{\lambda_2}, \ldots, r^{\lambda_r}\}$ so that each $i$ appears $\lambda_i$ times in $s_\lambda$. For example, if $\lambda = 4211$ then $s_\lambda = \{1^4, 2^2, 3^1, 4^1\}$.

For any $s \subset \mathbb{N}$ with $\lambda(s) = \lambda$, we reindex $s$ so that $s = \tilde{s} := \{\tilde{s}_1^{\lambda_1}, \ldots, \tilde{s}_r^{\lambda_r}\}$, with the convention $\tilde{s}_i < \tilde{s}_{i+1}$ whenever $\lambda_i = \lambda_{i+1}$. We define the canonical mapping $\sigma_{s,\lambda} : [r] \to s$ by $\sigma_{s,\lambda}(i) = \tilde{s}_i$ for each $i = 1, \ldots, r$. For example, let $s = \{1^2, 3^1, 4^2, 5^4\}$ so that $\lambda(s) = 4221$ and $s_\lambda = \{1^4, 2^2, 3^2, 4^1\}$. Then we reindex $s$ to obtain $\tilde{s} = \{5^4, 1^2, 4^2, 3^1\}$ and $\sigma_{s,\lambda} : [4] \to s$ assigns $\sigma_{s,\lambda}(1) = 5$, $\sigma_{s,\lambda}(2) = 1$, $\sigma_{s,\lambda}(3) = 4$, and $\sigma_{s,\lambda}(4) = 3$, yielding $s^{\sigma_{s,\lambda}} := \{\sigma_{s,\lambda}^{-1}(s_1)^{m_1}, \ldots, \sigma_{s,\lambda}^{-1}(s_r)^{m_r}\} = s_\lambda$.

For every $s \subset \mathbb{N}$, we define $\|\cdot\|_s$ by

$$\|M\|_s = (\|M^*_{s,1}\|, \ldots, \|M^*_{s,k}\|), \tag{28}$$

the limit of $M^*_s$ as an $\mathcal{L}^\lambda$-structure, where $\|M^*_{s,j}\|$ is the combinatorial limit of $M^*_{s,j}$ when encoded as an $(i_j - q)^{k_j}$-structure, with any prespecified convention for ordering the components of $M^*_{s,j} = (M^*_{s,1}, \ldots, M^*_{s,k_j})$. We write

$$\|M\|_s = \mathbf{0} \quad \text{if and only if} \quad \|M^*_{s,j}\| = \mathbf{0}_{(i_j-q)^{k_j}} \quad \text{for all } j = 1, \ldots, k, \tag{29}$$

where recall $\mathbf{0}_\mathcal{L}$ is the combinatorial limit of the empty structure $\mathbf{0}^\mathcal{L}_\mathbb{N}$ with signature $\mathcal{L}$, that is, $\mathbf{0}_\mathcal{L}$ is the probability measure which assigns unit mass to the empty structure $\mathbf{0}^\mathcal{L}_\mathbb{N}$. For $i_j < q$, we define $\|M^*_{s,j}\| = \mathbf{0}_{()}$, the limit of the null structure.

The above preparation anticipates Theorem 4.18 in which we decompose exchangeable $\sigma$-finite measures on $\mathcal{L}_\mathbb{N}$ according to how they treat various substructures. Below we write $\mu_\lambda$ to denote a measure on $\mathcal{L}_\mathbb{N}$ that satisfies (22),

$$M^*_{s_\lambda} = M \quad \text{for } \mu_\lambda\text{-almost every } M \in \mathcal{L}_\mathbb{N}, \tag{30}$$

and

$$s_M := \bigcap \{s' \subseteq \mathbb{N} : \|M\|_{s'} \neq \mathbf{0}\} = s_\lambda \quad \text{for } \mu_\lambda\text{-almost every } M \in \mathcal{L}_\mathbb{N}, \tag{31}$$

where the intersection of multisets is defined at the beginning of Section 4.3.

**Remark 4.17.** *Condition* (30) *says that $\mu_\lambda$ puts all of its support on the $s_\lambda$-component of $M \in \mathcal{L}_\mathbb{N}$. Condition* (31) *ensures that the support of $\mu_\lambda$ not on a proper subset of $s_\lambda$-indexed structures. Together, Conditions* (30) *and* (31) *are needed for identifiability purposes.*

We then define

$$\mu^*_\lambda(\cdot) = \sum_{s \subset \mathbb{N} : \lambda(s) = \lambda} \mu_\lambda(\{M \in \mathcal{L}_\mathbb{N} : M^{\sigma_{s,\lambda}^{-1}} \in \cdot\}) \tag{32}$$

for $\sigma_{s,\lambda}$ as defined above. For example, let $\lambda = 1$ be the only partition of integer 1 and put

$$\mu_\lambda((\mathbb{N}; \{i\})) = \begin{cases} c, & i = 1, \\ 0, & \text{otherwise}, \end{cases}$$

for some $c > 0$. Then $s_\lambda = \{1^1\}$ and $\mu_\lambda$ satisfies (22), (30), and (31). For any $k' > 1$, we note that $s = \{k'^1\}$ has $\lambda(s) = 1$ and $\sigma_{s,\lambda}(1) = k'$; whence, $M^{\sigma_{s,\lambda}^{-1}} = (\mathbb{N}; \{1\})$ if and only if $M = (\mathbb{N}; \{k'\})$ for some $k' \geq 1$. In this case, $\mu_\lambda^*$ assigns mass $c$ to each singleton subset $(\mathbb{N}; \{k'\})$, is exchangeable, and satisfies (22). Compare the definition of $\mu_1^*$ to that of $c \sum_{i=1}^\infty \epsilon_i$ in Theorem 2.7.

**Theorem 4.18** (Lévy–Itô–Khintchine representation for combinatorial Lévy processes). *Let $\mathcal{L} = (i_1, \ldots, i_k)$ be any signature and $\mathbf{X} = (X_t, t \geq 0)$ be an exchangeable combinatorial Lévy process on $\mathcal{L}_\mathbb{N}$ with rate measure $\mu$. Then there exists a unique measure $\nu_0$ on $\mathcal{E}_\mathcal{L}$ satisfying*

$$(33) \qquad \nu_0(\{\mathbf{0}_\mathcal{L}\}) = 0 \quad \text{and} \quad \int_{\mathcal{E}_\mathcal{L}} (1 - W(\{\mathbf{0}_{[i_k]}^\mathcal{L}\}))\nu_0(dW) < \infty,$$

*and measures $\mu_\lambda$ on $\mathcal{L}_\mathbb{N}$ satisfying (22), (30), and (31) for each $\lambda \vdash q$, $q = 1, \ldots, i_k$, such that*

$$(34) \qquad \mu = \nu_0^* + \sum_{q=1}^{i_k} \sum_{\lambda \vdash q} \mu_\lambda^*,$$

*where $\mu_\lambda^*$ is defined in (32).*

We call (34) the *Lévy–Itô–Khintchine representation* for exchangeable combinatorial Lévy processes. In a precise sense, see Theorem 4.19, $\nu_0^*$ describes the jump component while the $\mu_\lambda^*$ decompose the continuous component of $\|\mathbf{X}\|$.

**Theorem 4.19.** *Let $\mathbf{X} = (X_t, t \geq 0)$ be an exchangeable combinatorial Lévy process on $\mathcal{L}_\mathbb{N}$. Then the projection $\|\mathbf{X}\| = (\|X_t\|, t \geq 0)$ into $\mathcal{E}_\mathcal{L}$ exists almost surely and is a Feller process.*

The decomposition of the characteristic measure $\mu$ in (34) leads to a classification of the sample path behavior of the projected process $\|\mathbf{X}\| = (\|X_t\|, t \geq 0)$. We extend Definition 3.9 to processes $\mathbf{X} = (X_t, t \geq 0)$, calling $\mathbf{X}$ *dissociated* if $\mathbf{X}^S = (X_t|_S, t \geq 0)$ and $\mathbf{X}^{S'} = (X_t|_{S'}, t \geq 0)$ are independent for all $S, S' \subseteq \mathbb{N}$ for which $S \cap S' = \emptyset$.

**Theorem 4.20.** *Let $\mathbf{X} = (X_t, t \geq 0)$ be an exchangeable, dissociated combinatorial Lévy process on $\mathcal{L}_\mathbb{N}$. Then there is a deterministic, continuous path $t \mapsto Y_t$ such that $\|\mathbf{X}\| =_\mathcal{D} \mathbf{Y} = (Y_t, t \geq 0)$.*

**Corollary 4.21.** *Let $\mathbf{X} = (X_t, t \geq 0)$ be an exchangeable combinatorial Lévy process on $\mathcal{L}_\mathbb{N}$ with characteristic measure $\mu = \nu_0^* + \sum_{q=0}^{i_k} \sum_{\lambda \vdash q} \mu_\lambda^*$ as in (34). Then the sample paths of $\|\mathbf{X}\|$ are continuous except possibly at the times of jumps from the $\nu_0^*$ measure.*

Much of our remaining effort is dedicated to proving Theorems 4.18, 4.19, and 4.20, but first we illustrate the above theorems in specific, concrete cases.

## 5. Examples

**5.1. Set-valued Lévy processes.** In Section 2, we discussed combinatorial Lévy processes in the special case when $\mathcal{L} = (1)$ and $\mathbf{X} = (X_t, t \geq 0)$ evolves on the space of subsets of $\mathbb{N}$.

In this case, the combinatorial limit of $(\mathbb{N}; A)$ is determined by the limiting frequency of elements in a subset $A \subseteq \mathbb{N}$,
$$\pi(A) = \lim_{n \to \infty} n^{-1}|A \cap [n]|,$$
which is sure to exist with probability 1 by the strong law of large numbers and de Finetti's theorem. de Finetti's theorem also implies that the marginal distribution of $\mathbf{X}$ at any fixed time $t \geq 0$ is determined by a unique probability measure on $[0, 1]$ as in (5).

In the context of Theorem 4.18, the behavior of $\mathbf{X}$ on $\mathcal{L}_\mathbb{N}$ is described by a measure $\mu = \nu^* + \mathbf{c} \sum_{i=1}^\infty \epsilon_i$ with components defined as in Theorem 2.7. The first component $\nu^*$ is induced from a measure $\nu$ satisfying (6), the analog to (33) in the special case of set-valued processes. The second component $\mathbf{c} \sum_{i=1}^\infty \epsilon_i$ plays the role of $\mu_1^*$ in (34) since $\lambda = 1$ is the only partition of the integer 1. The only nontrivial measures $\mu_1$ that satisfy (22), (30), and (31) for $\lambda = 1$ must be of the form
$$\mu_1((\mathbb{N}; A)) = \begin{cases} \mathbf{c}, & A = \{1\}, \\ 0, & \text{otherwise}. \end{cases}$$

Our definition of $\mu_\lambda^*$ in (32) gives $\mu_1^*(\cdot) = \mathbf{c} \sum_{i=1}^\infty \epsilon_i(\cdot)$. The contribution of $\mu_1^*$ to the characteristic measure of $\mathbf{X}$ is as discussed previously: each $i \in \mathbb{N}$ changes status independently at rate $\mathbf{c} \geq 0$, while the rest of $\mathbf{X}$ remains unchanged.

For the effect on the limit process $\|\mathbf{X}\|$, we need only consider how the projection to the simplex $((\pi(X_t), 1 - \pi(X_t)), t \geq 0)$ behaves. In this case, suppose $\mu = \sum_{i \geq 1} \epsilon_i$ has no $\nu^*$ component, that is, $\nu \equiv 0$. If $(\pi(X_t), 1 - \pi(X_t)) = (1/2, 1/2)$ and $s > 0$, then each index has probability $1 - e^{-s}$ to change status at least once in the next $s$ units of time, independently of all the others. Some will change status two, three, four times, etc. during this period, but we assume $s$ is small enough to make those cases negligible. By independence and $(\pi(X_t), 1 - \pi(X_t)) = (1/2, 1/2)$, roughly half of the chosen elements will change from inside $X_t$ to outside $X_{t+s}$ and half will go from being outside $X_t$ to being inside $X_{t+s}$. The net effect is a wash, so that $(\pi(X_{t+s}), 1 - \pi(X_{t+s})) = (1/2, 1/2)$ for all $s > 0$ with probability 1. Assuming $X_0 = \emptyset$, the process begins at $(0, 1)$ and moves toward its steady state at $(1/2, 1/2)$ along the path $e^{-2t}$, that is, $(\pi(X_t), 1 - \pi(X_t)) = (\frac{1}{2}(1 - e^{-2t}), \frac{1}{2}(1 + e^{-2t}))$ for all $t > 0$ with probability 1.

To understand this more precisely, note that the infinitesimal jump rates governing the transitions of each element are given by the Q-matrix
$$Q = \begin{pmatrix} -1 & 1 \\ 1 & -1 \end{pmatrix}.$$

Thus, the marginal distribution of the state of each $i \in \mathbb{N}$ at time $t > 0$ is given by the matrix exponential
$$\begin{aligned}
e^{tQ} &= \sum_{k \geq 0} (tQ)^k/k! \\
&= \begin{pmatrix} 1 & 0 \\ 0 & 1 \end{pmatrix} + \frac{1}{2} \sum_{k \geq 1} \frac{t^k}{k!} \begin{pmatrix} (-1)^k 2^k & (-1)^{k+1} 2^k \\ (-1)^{k+1} 2^k & (-1)^k 2^k \end{pmatrix} \\
&= \begin{pmatrix} 1 & 0 \\ 0 & 1 \end{pmatrix} + \frac{1}{2} \begin{pmatrix} e^{-2t} - 1 & 1 - e^{-2t} \\ 1 - e^{-2t} & e^{-2t} - 1 \end{pmatrix} \\
&= \begin{pmatrix} \frac{1}{2}(1 + e^{-2t}) & \frac{1}{2}(1 - e^{-2t}) \\ \frac{1}{2}(1 - e^{-2t}) & \frac{1}{2}(1 + e^{-2t}) \end{pmatrix}.
\end{aligned}$$

It follows that

$$\mathbb{P}\{i \in X_t \mid X_0 = \emptyset\} = e^{tQ}(1, 2) = \frac{1}{2}(1 - e^{-2t}),$$

independently for each $i \in \mathbb{N}$. By the strong law of large numbers, we have $\pi(X_t) = \frac{1}{2}(1 - e^{-2t})$ a.s. for each $t \geq 0$.

I defer further discussion of the discontinuities in $\|\mathbf{X}\|$ in general combinatorial Lévy processes to the proof of Theorem 4.20 in Section 7.4.1.

5.2. **Graph-valued Lévy processes.** Let $\mathbf{X} = (X_t, t \geq 0)$ be a Lévy process on the space of directed graphs, possibly with self-loops, so that $\mathbf{X}$ evolves on the space of $\mathcal{L}$-structures with $\mathcal{L} = (2)$. By Theorem 4.18, the first component of $\mu$ in (34) is a measure $\nu_0$ on the space of graph limits satisfying (33). The second component decomposes according to the three partitions 1, 11, and 2 of the integers 1 and 2 as follows.

(1) $\mu_1$ is a measure on $\mathcal{L}_\mathbb{N}$ for which almost every $M = (\mathbb{N}, E)$ has $M = M^*_{\{1^1\}}$ and at least one of the conditions

$$\lim_{n \to \infty} n^{-1} \sum_{j=1}^n \mathbf{1}\{(1, j) \in E\} > 0 \quad \text{or} \quad \lim_{n \to \infty} n^{-1} \sum_{j=1}^n \mathbf{1}\{(j, 1) \in E\} > 0$$

holds.

(11) $\mu_{11}$ assigns 0 mass to all $M = (\mathbb{N}; E)$ except that for which at least one of $(1, 2) \in E$ and $(2, 1) \in E$ holds and $(i, j) \notin E$ otherwise.

(2) $\mu_2$ assigns 0 mass to all $M = (\mathbb{N}; E)$ except that for which $(1, 1) \in E$ and $(i, j) \notin E$ otherwise.

The jump rates of $\mathbf{X}$ are determined by $\mu = \nu^* + \mu_1^* + \mu_{11}^* + \mu_2^*$. At the time of a discontinuity in $\mathbf{X}$, either

(0) a strictly positive proportion of edges changes status according to a $\sigma$-finite measure $\nu_0^*$ on countable graphs,
(1) a positive proportion of edges incident to a specific vertex changes status and other edges stay fixed,
(11) edges involving a specific pair $\{i, j\}$, $i \neq j$, change status and the rest of the graph stays fixed, or
(2) a single self-loop $(i, i)$ changes status for a specific $i \in \mathbb{N}$ and the rest of the graph stays fixed.

In this special case, the limit process $\|\mathbf{X}\| = (\|X_t\|, t \geq 0)$ evolves on the space of *graph limits*. Lovász and Szegedy [34] introduced the term graph limit in 2006, but a more general concept originated with the Aldous–Hoover theorem [4, 28]; see [5, Theorem 14.11].

5.3. **Networks with a distinguished community.** Combining the structures in the previous two sections, we get signature $\mathcal{L} = (1, 2)$, which corresponds to a structure $M = (\mathbb{N}; A, E)$ with $A \subseteq \mathbb{N}$ and $E \subseteq \mathbb{N} \times \mathbb{N}$. In this case, a combinatorial Lévy process $\mathbf{X} = (X_t, t \geq 0)$ offers the interpretation as the evolution of a network along with a distinguished community of its vertices. As in the previous section, we must consider partitions of integers 1 and 2, so Theorem 4.18 characterizes exchangeable processes $\mathbf{X}$ by a $\sigma$-finite measure $\nu_0$ on $\mathcal{E}_\mathcal{L}$ and measures $\mu_1, \mu_{11}, \mu_2$. The $\nu_0$ measure governs a joint evolution of the community and the network such that atoms from $\nu_0$ cause a positive proportion of elements to change community status and/or a positive proportion of edges to change status. The $\mu_\lambda$ measures play a similar role to Section 5.2 with some modifications. For $\lambda = 1$, $\mu_1$ allows for the status of element 1 to change in the subset $A$ as well as a change to a positive proportion of edges incident to element 1 as in Section 5.2. For $\lambda = 11$, $\mu_{11}$ is just as in Section 5.2: there is a change to at least one of the edges $(1, 2)$ and $(2, 1)$ and no change in the community structure $A$. For $\lambda = 2$, $\mu_2$ allows for a change to the status of element 1 in the community structure as well as a change to the status of edge $(1, 1)$ in $E$.

## 6. Characterization of exchangeable $\sigma$-finite measures

### 6.1. Limits of combinatorial structures.
Recall definition (16) of the limiting densities of a structure $M$.

**Theorem 6.1** (Aldous–Hoover theorem for $\mathcal{L}$-structures [4, 28, 30]). *Let $\mathcal{L} = (i_1, \ldots, i_k)$ be a signature and $X$ be an exchangeable $\mathcal{L}$-structure over $\mathbb{N}$. Then there exists a measurable function $g = (g_1, \ldots, g_k)$ with $g_j : [0, 1]^{2^{i_j}} \to \{0, 1\}$ for each $j = 1, \ldots, k$ such that $X =_{\mathcal{D}} X^g = (\mathbb{N}; X_1^g, \ldots, X_k^g)$, where*

$$(35) \quad X_j^g(a) = g_j((\xi_s)_{s \leq a}), \quad a = (a_1, \ldots, a_{i_j}) \in \mathbb{N}^{i_j},$$

*for $(\xi_s)_{s \subset \mathbb{N} : |s| \leq i_k}$ a collection of i.i.d. Uniform$[0, 1]$ random variables. In particular, $X$ is conditionally dissociated given its tail $\sigma$-field.*

The key technical outcome of Theorem 6.1 is that every infinite exchangeable random structure is conditionally dissociated given its tail $\sigma$-field. (In the context of representation (35), the random variable $\xi_\emptyset$, which is common to all components, determines the tail $\sigma$-field. From this, we immediately see that, given $\xi_\emptyset$, the structure $X^g$ in (35) is conditionally dissociated.) The representation in (35) serves mostly as an analogy to our decomposition of exchangeable combinatorial Lévy processes in Theorem 4.18.

**Proposition 6.2.** *Let $X$ be an exchangeable $\mathcal{L}$-structure over $\mathbb{N}$. Then $(\delta(A; X), A \in \bigcup_{m \in \mathbb{N}} \mathcal{L}_{[m]})$ exists almost surely and determines a random probability measure $\|X\|$ on $\mathcal{L}_\mathbb{N}$.*

*Proof.* In addition to being exchangeable, we first assume that $X$ is dissociated, that is, $X|_S$ and $X|_T$ are independent whenever $S$ and $T$ are disjoint, as in Definition (3.9). For a fixed $\mathcal{L}$-structure $A = ([m]; A_1, \ldots, A_k)$ over $[m]$, we define

$$Z_n := \frac{1}{n^{\downarrow m}} \sum_{\varphi:[m] \to [n]} \mathbf{1}\{X|_{[n]}^\varphi = A\}, \quad \text{for each } n = 1, 2, \ldots.$$

For each $n \geq 1$, we define the $\sigma$-field $\mathcal{F}_n := \sigma\langle Z_{n+1}, Z_{n+2}, \ldots \rangle$. For any injection $\varphi : [m] \to [n + 1]$, exchangeability of $X$ implies

$$\mathbb{P}\{X|_{[n+1]}^\varphi = A \mid \mathcal{F}_n\} = \mathbb{P}\{X|_{[m]} = A \mid \mathcal{F}_n\};$$

whence,

$$\begin{aligned}
Z_{n+1} &= E\left( \frac{1}{(n+1)^{\downarrow m}} \sum_{\psi:[m] \to [n+1]} \mathbf{1}\{X|_{[n+1]}^\psi = A\} \mid \mathcal{F}_n \right) \\
&= \frac{1}{(n+1)^{\downarrow m}} \sum_{\psi:[m] \to [n+1]} E(\mathbf{1}\{X|_{[n+1]}^\psi = A\} \mid \mathcal{F}_n) \\
&= \mathbb{P}\{X|_{[m]} = A \mid \mathcal{F}_n\}.
\end{aligned}$$

Thus,

$$\begin{aligned} E(Z_n \mid \mathcal{F}_n) &= E\left( \frac{1}{n^{\downarrow m}} \sum_{\psi:[m]\to[n]} \mathbf{1}\{X|_{[n]}^{\psi} = A\} \mid \mathcal{F}_n \right) \\ &= \frac{1}{n^{\downarrow m}} \sum_{\psi:[m]\to[n+1]\text{ s.t. range}(\psi)\subseteq[n]} E(\mathbf{1}\{X|_{[n]}^{\psi} = A\} \mid \mathcal{F}_n) \\ &= \mathbb{P}(X|_{[m]} = A \mid \mathcal{F}_n) \\ &= Z_{n+1}, \end{aligned}$$

and $(Z_n, n \in \mathbb{N})$ is a reverse martingale with respect to the filtration $(\mathcal{F}_n, n \geq 1)$. By the reverse martingale convergence theorem, there exists a random variable $Z_\infty$ such that $Z_n \to Z_\infty$ almost surely. Since we have assumed $X$ is dissociated, the limit depends only on the tail $\sigma$-field $\mathcal{T} = \bigcap_{n\in\mathbb{N}} \mathcal{F}_n$ and, thus, is deterministic by the 0-1 law. That $\delta(A; X)$ exists for any exchangeable $X$ follows by the fact that any exchangeable $\mathcal{L}$-structure is conditionally dissociated given its tail $\sigma$-field, by Theorem 6.1. Almost sure existence of the infinite collection $(\delta(A; X), A \in \bigcup_{m\in\mathbb{N}} \mathcal{L}_{[m]})$ follows by countable additivity of probability measures.

To prove that $(\delta(A; X), A \in \bigcup_{m\in\mathbb{N}} \mathcal{L}_{[m]})$ determines a unique, exchangeable probability measure on $\mathcal{L}_\mathbb{N}$, we consider $A \in \mathcal{L}_{[m]}$ and $A' \in \mathcal{L}_{[n]}$ such that $A'|_{[m]} = A$, for $m \leq n$. For fixed $r \geq n$, the definition in (15) implies

$$\begin{aligned} \sum_{A'\in\mathcal{L}_{[n]}: A'|_{[m]}=A} \delta(A'; X|_{[r]}) &= \sum_{A'\in\mathcal{L}_{[n]}: A'|_{[m]}=A} \frac{1}{r^{\downarrow n}} \sum_{\varphi:[n]\to[r]} \mathbf{1}\{X|_{[r]}^{\varphi} = A'\} \\ &= \frac{1}{r^{\downarrow n}} \sum_{\varphi:[n]\to[r]} \sum_{A'\in\mathcal{L}_{[n]}: A'|_{[m]}=A} \mathbf{1}\{X|_{[r]}^{\varphi} = A'\} \\ &= \frac{1}{r^{\downarrow n}} \sum_{\varphi:[m]\to[r]} \mathbf{1}\{X|_{[r]}^{\varphi} = A\} \sum_{\text{extensions of } \varphi \text{ to } [n]\to[r]} 1 \\ &= \frac{1}{r^{\downarrow n}} \sum_{\varphi:[m]\to[r]} \mathbf{1}\{X|_{[r]}^{\varphi} = A\} \times (r-m)(r-m-1)\cdots(r-n+1) \\ &= \frac{1}{r^{\downarrow m}} \sum_{\varphi:[m]\to[r]} \mathbf{1}\{X|_{[r]}^{\varphi} = A\} \\ &= \delta(A; X|_{[r]}). \end{aligned}$$

Since $r \geq n$ is arbitrary, the probability measures induced on $\mathcal{L}_{[m]}$ and $\mathcal{L}_{[n]}$ are consistent for all $m \leq n$. Carathéodory's extension theorem implies an extension to a unique probability measure on $\mathcal{L}_\mathbb{N}$. Since each of the finite space distributions is exchangeable, so is the distribution induced on $\mathcal{L}_\mathbb{N}$. □

By Proposition 6.2, every infinite exchangeable $\mathcal{L}$-structure projects to a random limit in $\mathcal{E}_\mathcal{L}$. Conversely, by Theorem 6.1, the law of every infinite exchangeable $\mathcal{L}$-structure $X$ is determined by a probability measure $\nu$ on $\mathcal{E}_\mathcal{L}$ such that $X \sim \nu^*$, where $\nu^*$ is defined in (25). By the projective structure of $\mathcal{L}_\mathbb{N}$, $\nu^*$ is uniquely determined by the induced measures

$$\nu^{*(n)}(M) := \nu^*(\{M^* \in \mathcal{L}_\mathbb{N} : M^*|_{[n]} = M\}), \quad M \in \mathcal{L}_{[n]},$$

for every $n \in \mathbb{N}$.

6.2. **Characterization of measures.** We are especially interested in Lévy processes that evolve in continuous time on $\mathcal{L}_{\mathbb{N}}$ and, therefore, can jump possibly infinitely often in bounded time intervals. To see the additional possibilities in this case, let $\mathcal{L} = (1)$ so that $\mu$ is an exchangeable measure on subsets of $\mathbb{N}$. For $c > 0$, we define

$$\mu(\cdot) = c \sum_{i=1}^{\infty} \mathbf{1}\{(\mathbb{N}; \{i\}) \in \cdot\},$$

which assigns mass $c$ to the singleton subsets of $\mathbb{N}$ and, thus, has infinite total mass. For $n = 1, 2, \ldots$, the restriction of $\mu$ to $\mathcal{L}_{[n]}$ is

$$\mu^{(n)}(M) = \begin{cases} c, & M = ([n]; \{i\}) \text{ for } i = 1, \ldots, n, \\ \infty, & M = ([n]; \emptyset), \\ 0, & \text{otherwise}, \end{cases}$$

which is finite and exchangeable on $\mathcal{L}_{[n]} \setminus \{([n]; \emptyset)\}$. On the other hand, let $c, c' \geq 0$ and define

$$\mu(\cdot) = c \sum_{i=1}^{\infty} \mathbf{1}\{(\mathbb{N}; \{i\}) \in \cdot\} + c' \sum_{i=1}^{\infty} \sum_{j=i+1}^{\infty} \mathbf{1}\{(\mathbb{N}; \{i, j\}) \in \cdot\},$$

so that singletons have mass $c$ and doubletons have mass $c'$. For $n \in \mathbb{N}$,

$$\mu^{(n)}(M) = c \sum_{i=1}^{n} \mathbf{1}\{M = ([n]; \{i\})\} + c' \sum_{i=1}^{n} \sum_{j=n+1}^{\infty} \mathbf{1}\{M = ([n]; \{i\})\} + c' \sum_{i=1}^{n-1} \sum_{j=i+1}^{n} \mathbf{1}\{M = ([n]; \{i, j\})\},$$

which is finite on $\mathcal{L}_{[n]} \setminus \{\mathbf{0}_{[n]}^{\mathcal{L}}\}$ only if $c' = 0$. (The middle term in the above expression results because the restriction of any $(\mathbb{N}; \{i, n+j\})$ to $[n]$ is $([n]; \{i\})$, for every $j = n+1, n+2, \ldots$.) Immediately, $\mu$ satisfies (22) only if it assigns no mass to doubleton subsets. The same argument rules out tripletons, quadrupletons, and so on.

**Theorem 6.3.** Let $\mathcal{L} = (i_1, \ldots, i_k)$ be a signature and $\mu$ be an exchangeable measure on $\mathcal{L}_{\mathbb{N}}$ that satisfies (22). Then there exists a unique measure $\nu_0$ on $\mathcal{E}_{\mathcal{L}}$ satisfying (33) and measures $\mu_\lambda$ satisfying (22), (30), and (31) for each $\lambda \vdash q$, $q = 1, \ldots, i_k$, such that

(36) $$\mu = \nu_0^* + \sum_{q=1}^{i_k} \sum_{\lambda \vdash q} \mu_\lambda^*,$$

where $\mu_\lambda^*$ is defined in (32).

We first show that any $\mu$ constructed as in (36) satisfies (22).

**Proposition 6.4.** Let $\nu_0$ satisfy (33). Then $\nu_0^*$ defined in (25) satisfies (22).

*Proof.* The lefthand side of (22) follows immediately from the lefthand side of (33). For the righthand side of (22), we need to show

$$\nu_0^*(\{M \in \mathcal{L}_{\mathbb{N}} : M|_{[n]} \neq \mathbf{0}_{[n]}^{\mathcal{L}}\}) < \infty \quad \text{for all } n \in \mathbb{N}.$$

We note that

$$\{M \in \mathcal{L}_{\mathbb{N}} : M|_{[n]} \neq \mathbf{0}_{[n]}^{\mathcal{L}}\} = \bigcup_{s = \{s_1 < \cdots < s_{i_k}\} \subset [n]} \{M \in \mathcal{L}_{\mathbb{N}} : M|_s \neq \mathbf{0}_s^{\mathcal{L}}\},$$

because $M|_{[n]} = \mathbf{0}^{\mathcal{L}}_{[n]}$ only if $M|_s$ is empty for all $s \subset [n]$ with $|s| = i_k$. By exchangeability of $\nu^*_0$,

$$\nu^*_0(\{M \in \mathcal{L}_{\mathbb{N}} : M|_s \neq \mathbf{0}^{\mathcal{L}}_s\}) = \nu^*_0(\{M \in \mathcal{L}_{\mathbb{N}} : M|_{[i_k]} = \mathbf{0}^{\mathcal{L}}_{[i_k]}\}) = \int_{\mathcal{E}_{\mathcal{L}}} (1 - W(\{\mathbf{0}^{\mathcal{L}}_{[i_k]}\}))\nu_0(dW)$$

for every $s = \{s_1 < \cdots < s_{i_k}\} \subseteq [n]$. Thus,

$$\begin{aligned}
\nu^*_0(\{M \in \mathcal{L}_{\mathbb{N}} : M|_{[n]} \neq \mathbf{0}^{\mathcal{L}}_{[n]}\}) &= \nu^*_0\left(\bigcup_{s=\{s_1<\cdots<s_{i_k}\}\subset[n]} \{M \in \mathcal{L}_{\mathbb{N}} : M|_s \neq \mathbf{0}^{\mathcal{L}}_s\}\right) \\
&\leq \sum_{s=\{s_1<\cdots<s_{i_k}\}\subset[n]} \nu^*_0(\{M \in \mathcal{L}_{\mathbb{N}} : M|_s \neq \mathbf{0}^{\mathcal{L}}_s\}) \\
&\leq n^{i_k} \int_{\mathcal{E}_{\mathcal{L}}} (1 - W(\{\mathbf{0}^{\mathcal{L}}_{[i_k]}\}))\nu_0(dW) \\
&< \infty
\end{aligned}$$

by the righthand side of (33). The proof is complete. □

**Proposition 6.5.** *Let $\mathcal{L} = (i_1, \ldots, i_k)$ be a signature and suppose that $\mu_\lambda$ is a measure on $\mathcal{L}_{\mathbb{N}}$ satisfying (22), (30), and (31) for some $\lambda \vdash q$ with $q = 1, \ldots, i_k$. Then $\mu^*_\lambda$, as defined in (32), satisfies (22).*

*Proof.* Let $\lambda \vdash q$ for some $q = 1, \ldots, i_k$. By (30), $M = M^*_{s_\lambda}$ for $\mu_\lambda$-almost every $M \in \mathcal{L}_{\mathbb{N}}$; whence, $\mu_\lambda(\{M \in \mathcal{L}_{\mathbb{N}} : M^{\sigma^{-1}_{s,\lambda}}|_{[n]} \neq \mathbf{0}^{\mathcal{L}}_{[n]}\}) = 0$ for all $s = \{s_1, \ldots, s_r\} \not\subseteq [n]$. For $s = \{s_1^{m_1}, \ldots, s_r^{m_r}\}$ such that $\{s_1, \ldots, s_r\} \subset [n]$, we observe that $M^{\sigma^{-1}_{s,\lambda}}|_{[n]} = M|_{[n]}^{\sigma^{-1}_{s,\lambda}}$ and, therefore,

$$\mu_\lambda(\{M \in \mathcal{L}_{\mathbb{N}} : M^{\sigma^{-1}_{s,\lambda}}|_{[n]} \neq \mathbf{0}^{\mathcal{L}}_{[n]}\}) = \mu_\lambda(\{M \in \mathcal{L}_{\mathbb{N}} : M|_{[n]}^{\sigma^{-1}_{s,\lambda}} \neq \mathbf{0}^{\mathcal{L}}_{[n]}\}) = \mu_\lambda(\{M \in \mathcal{L}_{\mathbb{N}} : M|_{[n]} \neq \mathbf{0}^{\mathcal{L}}_{[n]}\}).$$

It follows that

$$\begin{aligned}
\mu^*_\lambda(\{M \in \mathcal{L}_{\mathbb{N}} : M|_{[n]} \neq \mathbf{0}^{\mathcal{L}}_{[n]}\}) &= \sum_{s\subset\mathbb{N}:\lambda(s)=\lambda} \mu_\lambda(\{M \in \mathcal{L}_{\mathbb{N}} : M^{\sigma^{-1}_{s,\lambda}}|_{[n]} \neq \mathbf{0}^{\mathcal{L}}_{[n]}\}) \\
&= \sum_{s\subseteq[n]:\lambda(s)=\lambda} \mu_\lambda(\{M \in \mathcal{L}_{\mathbb{N}} : M|_{[n]}^{\sigma^{-1}_{s,\lambda}} \neq \mathbf{0}^{\mathcal{L}}_{[n]}\}) \\
&= \sum_{s\subseteq[n]:\lambda(s)=\lambda} \mu_\lambda(\{M \in \mathcal{L}_{\mathbb{N}} : M|_{[n]} \neq \mathbf{0}^{\mathcal{L}}_{[n]}\}) \\
&\leq n^q \mu_\lambda(\{M \in \mathcal{L}_{\mathbb{N}} : M|_{[n]} \neq \mathbf{0}^{\mathcal{L}}_{[n]}\}) \\
&< \infty
\end{aligned}$$

which establishes the righthand side of (33). The lefthand side of (33) follows from assuming that $\mu_\lambda$ satisfies (22), the fact that there are at most countably many multisets $s$ for which $\lambda(s) = \lambda$, and countable additivity of measures. □

Theorem 6.3 states the converse of Propositions 6.4 and 6.5: the above construction is true of every exchangeable $\sigma$-finite measure on $\mathcal{L}_{\mathbb{N}}$. We prove Theorem 6.3 in several steps.

**Lemma 6.6.** *Let $\mu$ be an exchangeable $\sigma$-finite measure on $\mathcal{L}_{\mathbb{N}}$. Then $\|M\|_s$ exists for all $s = \{s_1^{m_1}, \ldots, s_r^{m_r}\} \subset \mathbb{N}$ with $0 \leq |s| \leq i_k$, for $\mu$-almost every $M \in \mathcal{L}_{\mathbb{N}}$.*

*Proof.* Recall Definition 3.8 of combinatorial limit. Fix $n = 1, 2, \ldots$ and define $\mu_n$ as the restriction of $\mu$ to the event $\{M \in \mathcal{L}_\mathbb{N} : M|_{[n]} \neq \mathbf{0}^\mathcal{L}_{[n]}\}$. By the righthand side of (22), $\mu_n$ is finite, because

$$\mu_n(S) = \mu(\{M \in \mathcal{L}_\mathbb{N} : M \in S \text{ and } M|_{[n]} \neq \mathbf{0}^\mathcal{L}_{[n]}\}) \leq \mu(\{M \in \mathcal{L}_\mathbb{N} : M|_{[n]} \neq \mathbf{0}^\mathcal{L}_{[n]}\}) < \infty$$

for all measurable sets $S \subseteq \mathcal{L}_\mathbb{N}$, by the righthand side of (22). Furthermore, exchangeability of $\mu$ implies that $\mu_n$ is invariant with respect to permutations of $\mathbb{N}$ that fix $[n]$. We define the shifted measure $\overleftarrow{\mu}_n$ as the image of $\mu_n$ by $M \mapsto \overleftarrow{M}^n := (\mathbb{N}; \overleftarrow{M}^n_1, \ldots, \overleftarrow{M}^n_k)$, where

$$(37) \qquad (a_1, \ldots, a_{i_j}) \in \overleftarrow{M}^n_j \quad \text{if and only if} \quad (a_1 + n, \ldots, a_{i_j} + n) \in M_j,$$

for each $j = 1, \ldots, k$. We call $\overleftarrow{M}^n$ the *$n$-shift of $M$*. For example, if $M = (\mathbb{N}; \{1, 2, 5, 6, 8\})$, then $\overleftarrow{M}^1 = (\mathbb{N}; \{1, 4, 5, 7\})$, $\overleftarrow{M}^2 = (\mathbb{N}; \{3, 4, 6\})$, and so on.

The $n$-shifted measure $\overleftarrow{\mu}_n$ is exchangeable and finite; therefore, $\overleftarrow{\mu}_n$ is proportional to an exchangeable probability measure and Proposition 6.2 implies that $\overleftarrow{\mu}_n$-almost every $M \in \mathcal{L}_\mathbb{N}$ possesses a unique combinatorial limit $\|M\|$. Furthermore, $\overleftarrow{\mu}_n$ induces a unique finite measure $\|\overleftarrow{\mu}_n\|$ on $\mathcal{E}_\mathcal{L}$ by

$$\|\overleftarrow{\mu}_n\|(\cdot) := \overleftarrow{\mu}_n(\{M \in \mathcal{L}_\mathbb{N} : \|M\| \in \cdot\}).$$

Now, suppose $\|\overleftarrow{M}^n\| = (\delta(A; \overleftarrow{M}^n), A \in \bigcup_{n \geq 1} \mathcal{L}_{[n]})$ exists for $n \geq 1$. Then for every $A \in \bigcup_{n \geq 1} \mathcal{L}_{[n]}$, we observe

$$\delta(A; \overleftarrow{M}^n) =$$

$$= \lim_{k \to \infty} \frac{1}{k^{\downarrow m}} \sum_{\varphi:[m] \to [k]} \mathbf{1}\{\overleftarrow{M}^n|^\varphi_{[k]} = A\}$$

$$= \lim_{k \to \infty} \frac{1}{k^{\downarrow m}} \sum_{\varphi:[m] \to [k+n]\setminus[n]} \mathbf{1}\{M|^\varphi_{[n+k]} = A\}$$

$$= \lim_{k \to \infty} \frac{(n+k)^{\downarrow m}}{k^{\downarrow m}} \frac{1}{(n+k)^{\downarrow m}} \sum_{\varphi:[m] \to [k+n]} \mathbf{1}\{M|^\varphi_{[k+n]} = A\} - \lim_{k \to \infty} \frac{1}{k^{\downarrow m}} \sum_{\varphi:[m] \to [n+k] \text{ s.t. } \varphi(j) \leq n \text{ some } j} \mathbf{1}\{M|^\varphi_{[n+k]} = A\}.$$

The first term above equals $\delta(A; M)$. The sandwich lemma shows that the second term converges to 0 by noticing that

$$0 \leq \frac{1}{k^{\downarrow m}} \sum_{\varphi:[m] \to [n+k] \text{ s.t. } \varphi(j) \leq n \text{ some } j} \mathbf{1}\{M|^\varphi_{[n+k]} = A\} \leq \frac{n(k+n)^{\downarrow(m-1)}}{k^{\downarrow m}}$$

for all $k \geq 1$. We conclude that $\delta(A; \overleftarrow{M}^n) = \delta(A; M)$ for all $A \in \bigcup_{n \geq 1} \mathcal{L}_{[n]}$ and, thus, the combinatorial limit $\|M\|$ depends only on the $n$-shift $\overleftarrow{M}^n$ for every $n \in \mathbb{N}$. It follows that $\mu_n$-almost every $M \in \mathcal{L}_\mathbb{N}$ possesses a limit as well. Finally, notice that the events $\{M \in \mathcal{L}_\mathbb{N} : M|_{[n]} \neq \mathbf{0}^\mathcal{L}_{[n]}\}$ increase to $\{M \in \mathcal{L}_\mathbb{N} : M \neq \mathbf{0}^\mathcal{L}_\mathbb{N}\}$ as $n \to \infty$. Since we have assumed that $\mu$ assigns no mass to $\{\mathbf{0}^\mathcal{L}_\mathbb{N}\}$, the monotone convergence theorem implies that $\mu_n \uparrow \mu$ as $n \to \infty$, and thus $\mu$-almost every $M \in \mathcal{L}_\mathbb{N}$ possesses a limit $\|M\|$.

The above argument shows that $\|M\|_s$ exists for $\mu$-almost every $M$ when $s = \emptyset$. The argument is similar for $s = \{s_1^{m_1}, \ldots, s_r^{m_r}\} \subset \mathbb{N}$ with $|s| = 1, \ldots, i_k$. In this case, we define $\mu_{n,s}$ to be the measure induced by $\mu_n$ through the map $M \mapsto M_s^*$. For any measurable set $S \subseteq \mathcal{L}_\mathbb{N}$, we define $S_s^* := \{M_s^* : M \in S\}$, and we see that

$$\mu_{n,s}(S_s^*) = \mu(\{M \in \mathcal{L}_\mathbb{N} : M_s^* \in S_s^* \text{ and } M|_{[n]} \neq \mathbf{0}^{\mathcal{L}}_{[n]}\}) \leq \mu(\{M \in \mathcal{L}_\mathbb{N} : M|_{[n]} \neq \mathbf{0}^{\mathcal{L}}_{[n]}\}) < \infty.$$

Furthermore, $\mu_{n,s}$ is invariant with respect to permutations that fix $[n]$ and $\{s\}$. Taking $n \geq 1 + \max_{1 \leq j \leq r} s_r$, we may define the $n$-shift measure $\overleftarrow{\mu}_{n,s}$ just as before, so that $\overleftarrow{\mu}_{n,s'}$ is exchangeable and finite for every $s' \subseteq [n]$. It follows that $\|M\|_s$ exists for $\overleftarrow{\mu}_{n,s}$-almost every $M \in \mathcal{L}_\mathbb{N}$. Since we have defined $\|M\|_s$ to depend only on $M_s^*$, it follows that $\|M\|_s$ exists for $\overleftarrow{\mu}_n$-almost every $M \in \mathcal{L}_\mathbb{N}$ and, by monotone convergence, $\mu$-almost every $M \in \mathcal{L}_\mathbb{N}$. □

When $\mu$ is a probability measure, the meaning of the limit $\|M\|$ for $M \sim \mu$ is determined by our discussion surrounding Definition 3.8. In this case, $\|M\|$ describes the limiting homomorphism densities of $M$ and describes an alternative way of generating $M \sim \mu$ in accordance with the Aldous–Hoover theorem: first generate $\|M\|$ and, given $\|M\| = W$, sample $M$ from the random probability measure $W$. This is exactly the description in (25).

However, our running example for subsets shows why this interpretation is not available for $\sigma$-finite measures $\mu$; see the discussion in Section 5. In particular, $\mu$ can assign positive mass to singleton subsets $M = (\mathbb{N}; \{i\})$, for every $i = 1, 2, \ldots$, in which case $\|M\| = \mathbf{0}_{\mathcal{L}}$ but $M \neq (\mathbb{N}; \emptyset) = \mathbf{0}^{\mathcal{L}}_\mathbb{N}$. To characterize $\mu$, we must understand how it treats events of the form $\{M \in \mathcal{L}_\mathbb{N} : \|M\| = \mathbf{0}_{\mathcal{L}}\}$.

**Lemma 6.7.** *Suppose $\mu$ is exchangeable and $\sigma$-finite, and suppose $\mu$-almost every $M \in \mathcal{L}_\mathbb{N}$ has $\|M\| \neq \mathbf{0}_{\mathcal{L}}$. Then there exists a unique measure $\nu$ on $\mathcal{E}_{\mathcal{L}}$ satisfying (33) such that $\mu = \nu^*$.*

*Proof.* As in Lemma 6.6, we let $\mu_n$ denote the restriction of $\mu$ to $\{M \in \mathcal{L}_\mathbb{N} : M|_{[n]} \neq \mathbf{0}^{\mathcal{L}}_{[n]}\}$ and we write $\overleftarrow{\mu}_n$ as the image of $\mu_n$ by the $n$-shift operation (37). Since $\overleftarrow{\mu}_n$ is exchangeable, the combinatorial limit $\|M\|$ exists for $\overleftarrow{\mu}_n$-almost every $M \in \mathcal{L}_\mathbb{N}$, allowing us to write

$$\overleftarrow{\mu}_n(\cdot) = \int_{\mathcal{E}_{\mathcal{L}}} W(\cdot) \|\overleftarrow{\mu}_n\|(dW).$$

By assumption, $\mu_n$-almost every $M$ has $\|M\| \neq \mathbf{0}_{\mathcal{L}}$, from which it follows that

$$\mu_n(\{\overleftarrow{\tilde{M}}^n|_{[i_k]} \neq \mathbf{0}^{\mathcal{L}}_{[i_k]}\}) = \int_{\mathcal{E}_{\mathcal{L}}} (1 - W(\{\mathbf{0}^{\mathcal{L}}_{[i_k]}\})) \|\overleftarrow{\mu}_n\|(dW),$$

where $\overleftarrow{\tilde{M}}^n$ is the $n$-shift from (37). Again, $\mu_n \uparrow \mu$ implies $\|\overleftarrow{\mu}_n\| \uparrow \nu$ for some measure $\nu$ on $\mathcal{E}_{\mathcal{L}}$; whence, $\nu(\{\mathbf{0}_{\mathcal{L}}\}) = 0$ and

$$\mu_n(\{\overleftarrow{\tilde{M}}^n|_{[i_k]} \neq \mathbf{0}^{\mathcal{L}}_{[i_k]}\}) \uparrow \int_{\mathcal{E}_{\mathcal{L}}} (1 - W(\{\mathbf{0}^{\mathcal{L}}_{[i_k]}\})) \nu(dW).$$

Furthermore, $\mu_n \uparrow \mu$ implies

$$\mu_n(\{\overleftarrow{\tilde{M}}^n|_{[i_k]} \neq \mathbf{0}^{\mathcal{L}}_{[i_k]}\}) \leq \mu(\{\overleftarrow{\tilde{M}}^n|_{[i_k]} \neq \mathbf{0}^{\mathcal{L}}_{[i_k]}\}) = \mu(\{M|_{[i_k]} \neq \mathbf{0}^{\mathcal{L}}_{[i_k]}\}) < \infty,$$

by the righthand side of (22) and exchangeability. It follows that

$$\int_{\mathcal{E}_{\mathcal{L}}} (1 - W(\{\mathbf{0}^{\mathcal{L}}_{[i_k]}\})) \nu(dW) < \infty,$$

so that $\nu$ satisfies (33).

To establish the identity $\mu = \nu^*$, we observe that

$$\mu(\{M^* \in \mathcal{L}_\mathbb{N} : M^*|_{[n]} = M, \|M^*\| \neq \mathbf{0}_\mathcal{L}\}) =$$
$$= \lim_{m \uparrow \infty} \mu(\{M^* \in \mathcal{L}_\mathbb{N} : M^*|_{[n]} = M, \|M^*\| \neq \mathbf{0}_\mathcal{L}, \overleftarrow{M}^{*n}|_{[m]} \neq \mathbf{0}^\mathcal{L}_{[m]}\}),$$

for every fixed $n \in \mathbb{N}$ and $M \in \mathcal{L}_{[n]} \setminus \{\mathbf{0}^\mathcal{L}_{[n]}\}$. The identity follows by the monotone convergence theorem because

$$\{M^* \in \mathcal{L}_\mathbb{N} : M^*|_{[n]} = M, \|M^*\| \neq \mathbf{0}_\mathcal{L}, \overleftarrow{M}^{*n}|_{[m]} \neq \mathbf{0}^\mathcal{L}_{[m]}\}$$

increases to

$$\{M^* \in \mathcal{L}_\mathbb{N} : M^*|_{[n]} = M, \|M^*\| \neq \mathbf{0}_\mathcal{L}, \overleftarrow{M}^{*n} \neq \mathbf{0}^\mathcal{L}_\mathbb{N}\} \quad \text{as } m \to \infty,$$

and $\{M^* \in \mathcal{L}_\mathbb{N} : \|M^*\| \neq \mathbf{0}_\mathcal{L}\} \subseteq \{M^* \in \mathcal{L}_\mathbb{N} : \overleftarrow{M}^{*n} \neq \mathbf{0}^\mathcal{L}_\mathbb{N}\}$ for every $n \in \mathbb{N}$. Exchangeability of $\mu$ allows us to permute the blocks $\{1, \ldots, n\}$ and $\{n+1, \ldots, n+m\}$ so that

$$\mu(\{M^* \in \mathcal{L}_\mathbb{N} : M^*|_{[n]} = M, \overleftarrow{M}^{*n}|_{[m]} \neq \mathbf{0}^\mathcal{L}_{[m]}, \|M^*\| \neq \mathbf{0}_\mathcal{L}\}) = \overleftarrow{\mu}_m(\{M^* \in \mathcal{L}_\mathbb{N} : M^*|_{[n]} = M, \|M^*\| \neq \mathbf{0}_\mathcal{L}\}).$$

Now, $\overleftarrow{\mu}_m$ is exchangeable and previous arguments imply

$$\overleftarrow{\mu}_m(\{M^* \in \mathcal{L}_\mathbb{N} : M^*|_{[n]} = M, \|M^*\| \neq \mathbf{0}_\mathcal{L}\}) = \int_{\mathcal{E}_\mathcal{L}} W(\{M^*|_{[n]} = M\}) \|\overleftarrow{\mu}_m\|(dW),$$

which converges to

$$\int_{\mathcal{E}_\mathcal{L}} W(\{M^*|_{[n]} = M\}) \nu(dW) = \nu^*(\{M^*|_{[n]} = M\}).$$

Since we chose $n$ arbitrarily and we restricted $\|M^*\| \neq \mathbf{0}_\mathcal{L}$ so that $M^* \neq \mathbf{0}^\mathcal{L}_\mathbb{N}$, we must have equality of $\mu$ and $\nu^*$ on the $\pi$-system that generates the Borel $\sigma$-field. Since $\sigma$-finite measures are determined by their behavior on a generating $\pi$-system, the proof is complete. □

**Lemma 6.8.** *Let $\mu$ be an exchangeable, $\sigma$-finite measure on $\mathcal{L}_\mathbb{N}$ for which $\mu$-almost every $M \in \mathcal{L}_\mathbb{N}$ has $\|M\| = \mathbf{0}_\mathcal{L}$. Then there are unique measures $\mu_\lambda$ satisfying (22), (30), and (31) for each $\lambda \vdash q$, $q = 1, \ldots, i_k$, such that*

$$\mu = \sum_{q=1}^{i_k} \sum_{\lambda \vdash q} \mu_\lambda^*. \tag{38}$$

The following proof appeals repeatedly to the exchangeability and $\sigma$-finiteness properties of $\mu$. The main idea is that any $\mu_\lambda$ satisfying (22), (30), and (31) is invariant under relabeling by all permutations $\sigma : \mathbb{N} \to \mathbb{N}$ that behave as the identity on $s_\lambda$. Consequently, if $\mu_\lambda$ assigns positive mass to the event $M^*_{s_\lambda} \neq M$ for some $M \in \mathcal{L}_\mathbb{N}$, then either $M$ does not satisfy $s_M = s_\lambda$ or $\mu_\lambda$ is not $\sigma$-finite. The proof proceeds by ruling out all possibilities for which one of these conditions must be violated.

*Proof.* Throughout we fix a signature $\mathcal{L} = (i_1, \ldots, i_k)$ with $i_k \geq 1$ and let $\mu$ be any $\sigma$-finite measure on $\mathcal{L}_\mathbb{N}$ assigning all its mass to the event $\|M\| = \mathbf{0}_\mathcal{L}$.

Fix any $M \in \mathcal{L}_\mathbb{N}$ and, for any multiset $s' \subset \mathbb{N}$, recall the definition of $\|M\|_{s'}$ in (28) and our definition of the event $\{\|M\|_{s'} = \mathbf{0}\}$ in (29). If $\|M\|_{s'}$ exists for all $s' \subset \mathbb{N}$, then we define $s_M = \bigcap \{s' \subset \mathbb{N} : \|M\|_{s'} \neq \mathbf{0}\}$ and $\lambda_M = \lambda(s_M)$, where $\lambda(s)$ is the partition induced by $s$ as in Section 4.3. By Lemma 6.6, $\|M\|_{s'}$ exists for $\mu$-almost every $M \in \mathcal{L}_\mathbb{N}$ whenever $\mu$ is $\sigma$-finite

and exchangeable; thus, $s_M$ and $\lambda_M$ are well defined $\mu$-almost everywhere on $\mathcal{L}_{\mathbb{N}}$. From now on, we tacitly assume that $s_M$ is well defined whenever we speak of a generic $M \in \mathcal{L}_{\mathbb{N}}$.

We propose the following candidates for the decomposition measures $\mu_\lambda$ in (38). Recall the definition of $s_\lambda$, where $\lambda \vdash q$ for some $q = 1, \ldots, i_k$, from Section 4.3. For every $\lambda \vdash q$, $q = 1, \ldots, i_k$, we define

$$\varrho_\lambda = \mu \mathbf{1}_{\{M \in \mathcal{L}_{\mathbb{N}} : s_M = s_\lambda\}} \tag{39}$$

to be the restriction of $\mu$ to structures $M$ for which $s_M = s_\lambda$. We define $\varrho_\lambda^*$ from $\varrho_\lambda$ as in (32) and we claim $\mu = \sum_{q=1}^{i_k} \sum_{\lambda \vdash q} \varrho_\lambda^*$.

First, by our definition of $\|M\|_{s'}$ in (28), it is clear that $\|M\|_{s'} = \mathbf{0}$ whenever $|s'| > i_k$; therefore, $|s_M| \leq i_k$ for all $M \in \mathcal{L}_{\mathbb{N}}$ and it is sufficient to decompose $\mu$ in terms of partitions of integers less than or equal to $i_k$. Second, we must rule out the possibility that $s_M = \emptyset$, since we have not defined $\varrho_\lambda$ for the empty partition of the integer 0. For this we define $\Sigma(M) := \{s' \subset \mathbb{N} : \|M\|_{s'} \neq \mathbf{0}\}$ for every $M \in \mathcal{L}_{\mathbb{N}}$, and so $s_M = \bigcap_{s' \in \Sigma(M)} s'$.

Suppose $\mu$ assigns positive measure to some $\mathcal{B} \subseteq \mathcal{L}_{\mathbb{N}}$ such that $s_M = \emptyset$ for all $M \in \mathcal{B}$. Then for every $M \in \mathcal{B}$ there exists $s', s'' \in \Sigma(M)$ such that $\{s'\} \cap \{s''\} = \emptyset$. We first note that if there are $s', s'' \in \Sigma(M)$ such that $|s'| = |s''| = q$ and $\{s'\} \cap \{s''\} = \emptyset$, then there must be a smallest $q = 1, \ldots, i_k$ for which this holds. We let $\mu_{q,n}$ be the restriction of $\mu$ to the event $E_q \cap \{M|_{[n]} \neq \mathbf{0}_{[n]}^{\mathcal{L}}\}$, where

$$E_q = \{M \in \mathcal{L}_{\mathbb{N}} : s(M) = \emptyset \text{ and } q \text{ is smallest s.t. } s', s'' \in \Sigma(M), |s'| = |s''| = q, \text{ and } \{s'\} \cap \{s''\} = \emptyset\}.$$

For any $M \in E_q$, we define

$$F(M) = \limsup_{n \to \infty} \frac{1}{n^q} \sum_{s' \subseteq [n] : |s'|=q} \mathbf{1}\{s' \in \Sigma(M)\}$$

to be the limit superior of the fraction of $s' \subseteq \mathbb{N}$ with $|s'| = q$ that are in $\Sigma(M)$. Since there are finitely many partitions of $q$, we have

$$F(M) \leq \sum_{\lambda \vdash q} \limsup_{n \to \infty} \frac{1}{n^q} \sum_{s' \subseteq [n] : \lambda(s')=\lambda} \mathbf{1}\{s' \in \Sigma(M)\},$$

so that if $F(M) > 0$, then there is at least one $\lambda \vdash q$ for which the limit superior of the fraction of $s' \subseteq \mathbb{N}$ with $\lambda(s') = \lambda$ and $s' \in \Sigma(M)$. If $s' \in \Sigma(M)$, then $\|M\|_{s'} \neq \mathbf{0}$, which implies that the $\|M\|_{s'}(\mathbf{0}_{\mathbb{N}}^{\lambda(s')}) < 1$, where here we regard $M_{s'}^*$ as an $\mathcal{L}^{\lambda(s')}$-structure and treat $\|M\|_{s'}$ as the combinatorial limit of $M_{s'}^*$; refer back to Section 4.3 and the discussion surrounding (28) and (29) for a detailed explanation of this notation.

Thus, if

$$\limsup_{n \to \infty} \frac{1}{n^q} \sum_{s' \subseteq [n] : \lambda(s')=\lambda} \mathbf{1}\{s' \in \Sigma(M)\} > 0$$

for some $\lambda \vdash q$, then $\|M\|(\mathbf{0}_{\mathbb{N}}^{\mathcal{L}}) < 1$, contradicting the assumption that $\mu$-almost every $M$ has $\|M\| = \mathbf{0}_{\mathcal{L}}$. It follows that the limit superior of the fraction of pairwise disjoint $s' \in \Sigma(M)$ with $|s'| = q$ must be 0 and, therefore, the limiting fraction exists and equals 0.

Suppose for now that there are exactly two $s', s'' \in \Sigma(M)$ for which $|s'| = |s''| = q$ and $\{s'\} \cap \{s''\} = \emptyset$. By exchangeability of $\mu$, if $M \in \mathcal{B}$ then $M^\sigma \in \mathcal{B}$ for all $\sigma : \mathbb{N} \to \mathbb{N}$. Under any such relabeling, the sets $s', s''$ are permuted accordingly, and so we may assume that $\max s' < \min s''$ without loss of generality. Furthermore, we can choose $n$ so that $\max s' \leq n < \min s''$, in which case the image of $\mu_{q,n}$ by the $n$-shift $\overleftarrow{\mu}_{q,n}$ assigns positive mass to $M \in \mathcal{L}_\mathbb{N}$ with a single $s'' \in \Sigma(M)$ having $|s''| = q$. Again, by exchangeability, $\overleftarrow{\mu}_{q,n}$ assigns positive measure to all $M \in \mathcal{L}_\mathbb{N}$ with a single $s^* \in \Sigma(M)$ for which $|s^*| = q$ and $\lambda(s^*) = \lambda(s'')$. As there are infinitely many such $s^*$, $\overleftarrow{\mu}_{q,n}$ must assign zero mass to all such events or else $\overleftarrow{\mu}_{q,n}$ would have infinite total mass, a contradiction. It follows that $\mu_{q,n}$ assigns zero mass to the set of $M \in E_q$ with exactly two nonoverlapping $s', s'' \in \Sigma(M)$ satisfying $|s'| = |s''| = q$. The same argument carries through for $M \in E_q$ for which there are more than two (but a zero limiting fraction of) mutually disjoint $s^*$ satisfying the condition, because in such a case we can choose any two $s', s'' \in \Sigma(M)$ and apply the above argument to achieve a contradiction.

Now, suppose $s', s'' \in \Sigma(M)$ satisfy $s' \cap s'' \neq \emptyset$ and $\{s'\} \neq \{s''\}$, so that $\#(s' \cap s''^c) \in [q-1]$. Then once again we can assume without loss of generality that $s' \cap s''^c \subset \mathbb{N} \setminus [n]$ so that $\overleftarrow{\mu}_{q,n}$ assigns positive measure to the event $\{M \in \mathcal{L}_\mathbb{N} : \|M\|_{s' \cap s''^c} \neq \mathbf{0}\}$. But again, exchangeability implies that $\overleftarrow{\mu}_{q,n}$ assigns positive measure to the event $\|M\|_{s^*}$ for all $s^*$ with $\lambda(s^*) = \lambda(s' \cap s''^c)$. As there are infinitely many such $s^*$ and $\overleftarrow{\mu}_{q,n}$ is finite, all such events must have measure 0 under $\overleftarrow{\mu}_{q,n}$ and, hence, also $\mu_{q,n}$. We conclude that for $\mu$-almost every $M$, there is a unique (possibly empty) subset $A_q$ for every $q = 1, \ldots, i_k$, such that $\{s'\} = A_q$ for all $s' \in \Sigma(M)$ having $|s'| = q$.

The above argument establishes that for $\mu$-almost every $M$, any $s', s'' \in \Sigma(M)$ with $|s'| = |s''|$ must also have $\{s'\} = \{s''\}$. Thus, if $\mu$ assigns positive mass to the event $s(M) = \emptyset$, there must be $s', s'' \in \Sigma(M)$ with $|s'| < |s''|$ such that $s' \cap s'' = \emptyset$. We can rule this out by a similar argument to the case $|s'| = |s''|$ above. In particular, we can assume that $\max s' \leq n < \min s''$ so that the image of $\mu$ under the $n$-shift assigns positive measure to $M'$ with $\|M'\|_{s''} \neq \mathbf{0}$. By exchangeability, the $n$-shift also assigns positive measure to $M'$ with $\|M'\|_{s^*} \neq \mathbf{0}$ for all $s^*$ having $\lambda(s^*) = \lambda(s'')$. Since there are infinitely many such $s^*$, finiteness of the $n$-shift measure once again forces each of these to have measure 0. We conclude that $\mu$-almost every $M \in \mathcal{L}_\mathbb{N}$ has $s(M) \neq \emptyset$ and, moreover, $\{s'\} = \{s''\}$ for all $s', s'' \in \Sigma(M)$.

It follows immediately that $\mu$ decomposes as

$$\mu = \sum_{q=1}^{i_k} \sum_{\lambda \vdash q} \mu \mathbf{1}_{\{M \in \mathcal{L}_\mathbb{N} : \lambda_M = \lambda\}},$$

since the $\mu \mathbf{1}_{\{M \in \mathcal{L}_\mathbb{N} : \lambda_M = \lambda\}}$ are mutually singular for different $\lambda$. We need to show that each $\mu \mathbf{1}_{\{M \in \mathcal{L}_\mathbb{N} : \lambda_M = \lambda\}}$ coincides with $\varrho_\lambda^*$ for $\varrho_\lambda$ defined in (39).

We first note that each $\varrho_\lambda$ automatically satisfies (31) by the definition $\varrho_\lambda := \mu \mathbf{1}_{\{M \in \mathcal{L}_\mathbb{N} : s_M = s_\lambda\}}$. We must show that $M_{s_\lambda}^* = M$ for $\varrho_\lambda$-almost every $M \in \mathcal{L}_\mathbb{N}$. Suppose $\varrho_\lambda$ assigns positive measure to the event $\{M_{s_\lambda}^* \neq M\}$. Referring to the definition of $M_s^*$ in (27), there must be some $j = 1, \ldots, k$ and some $a \in \mathbb{N}^{i_j}$ not satisfying either of the top two conditions with respect to $s_\lambda$ but for which $M_{s,j}^*(a) \neq 0$. We rule out all possibilities as follows.

(1) Suppose $|a| \leq |s_\lambda|$ and $a \leq s_\lambda$ but $\{a\} \subsetneq \{s_\lambda\}$. In this case $M_{a,j}^*$ can be represented as a $(0)^{k_j}$-structure that automatically has $M_{a,j}^*(a) \neq 0$. Thus, $\|M\|_a \neq \mathbf{0}$, $a \in \Sigma(M)$, and $s_M \subseteq a \subsetneq s_\lambda$, a contradiction.
(2) Suppose $|a| \leq |s_\lambda|$ and $a \not\leq s_\lambda$. Then again $M_{a,j}^*$ is a non-empty $(0)^{k_j}$-structure for which $\|M\|_a \neq \mathbf{0}$, implying $s_M \subseteq a \cap s_\lambda \subsetneq s_\lambda$, a contradiction.

(3) Suppose $|a| > |s_\lambda|$ and $s_\lambda \not\leq a$. By definition $\varrho_\lambda$ is invariant with respect to relabeling by permutations that act as the identity on $s_\lambda$. If $\{s_\lambda\} \subsetneq \{a\}$, then $a \in \Sigma(M)$ contains an element not in $s_\lambda$, contradicting our analysis above, which shows that $\{s'\} = \{s''\}$ for all $s', s'' \in \Sigma(M)$. If $\{a\} \subseteq \{s_\lambda\}$, then invariance of $\varrho_\lambda$ with respect to permutations that fix $s_\lambda$ would imply $\|M\|_a \neq \mathbf{0}$ and $s_M \subseteq a \cap s_\lambda \subsetneq s_\lambda$, a contradiction.

It follows that each $\varrho_\lambda$ defined in (39) satisfies (30) and (31). Finally, since $\lambda(s_M) = \lambda(s_{M^\sigma})$ for all $M \in \mathcal{L}_\mathbb{N}$ and permutations $\sigma : \mathbb{N} \to \mathbb{N}$, it follows that $\varrho_\lambda^* = \mu \mathbf{1}_{\{M \in \mathcal{L}_\mathbb{N} : \lambda_M = \lambda\}}$, proving that $\varrho_\lambda$ is uniquely determined by $\mu$. Finally, each $\varrho_\lambda^*$ inherits $\sigma$-finiteness and exchangeability from $\mu$ because $\lambda_M = \lambda_{M^\sigma}$ for all permutations $\sigma : \mathbb{N} \to \mathbb{N}$. The proof is complete. □

*Proof of Theorem 6.3.* This is a consequence of Proposition 6.4 and Lemmas 6.6, 6.7, and 6.8. □

## 7. Proofs of main theorems

Theorem 6.3 is key to our conclusions about infinite exchangeable combinatorial Lévy processes. In this section, we prove the main theorems from Section 4.

7.1. **Discrete time combinatorial Lévy processes.** Theorem 4.5 is immediate from Definition 4.2. We now prove Theorem 4.15, which deals with discrete time combinatorial Lévy processes that are exchangeable.

*Proof of Theorem 4.15.* Let $\mathbf{X} = (X_m, m \geq 0)$ be an exchangeable combinatorial Lévy process in discrete time. By definition, the increments process $(\Delta_m, m \geq 1)$ defined by $\Delta_m := X_m \triangle X_{m-1}$ is a sequence of independent, identically distributed structures. The increment operator $\triangle$ satisfies

$$\Delta_m^\sigma = (X_m \triangle X_{m-1})^\sigma = X_m^\sigma \triangle X_{m-1}^\sigma. \tag{40}$$

By exchangeability of $\mathbf{X}$, we observe $\Delta_m^\sigma =_\mathcal{D} \Delta_m$ for all permutations $\sigma : \mathbb{N} \to \mathbb{N}$, from which exchangeability of the increments follows.

*Proof of Theorem 4.15: Infinite case.* Suppose $\mathbf{X}$ evolves on $\mathcal{L}_\mathbb{N}$ for some signature $\mathcal{L} = (i_1, \ldots, i_k)$ with $i_k \geq 1$. The representation by a unique probability measure $\nu$ on $\mathcal{E}_\mathcal{L}$ follows from Proposition 6.2 and exchangeability of the increments in (40).

Almost sure existence of $\|\mathbf{X}\|$ follows from Proposition 6.2 and countable additivity of probability measures. To establish the Markov property for $\|\mathbf{X}\|$, we introduce another process $\mathbf{P} = (P_m, m \geq 0)$ which couples with both $\mathbf{X}$ and $\|\mathbf{X}\|$ as follows.

Every $\mu \in \mathcal{E}_\mathcal{L}$ determines a unique transition probability measure $P_\mu$ on $\mathcal{L}_\mathbb{N}$ by

$$P_\mu(M, \cdot) := \mu(\{\Delta \in \mathcal{L}_\mathbb{N} : M \triangle \Delta \in \cdot\}), \quad M \in \mathcal{L}_\mathbb{N}. \tag{41}$$

This transition measure is exchangeable in the sense that $P_\mu(M^\sigma, \mathcal{B}^\sigma) = P_\mu(M, \mathcal{B})$ for all $M \in \mathcal{L}_\mathbb{N}, \mathcal{B} \subseteq \mathcal{L}_\mathbb{N}$, and permutations $\sigma : \mathbb{N} \to \mathbb{N}$. We write $\mathfrak{S}_\mathcal{L}$ to denote the set of all such transition probabilities on $\mathcal{L}_\mathbb{N}$.

Every $P \in \mathfrak{S}_\mathcal{L}$ acts on $\mathfrak{S}_\mathcal{L}$ by $P' \mapsto P \circ P'$, where

$$(P \circ P')(M, \cdot) := \int_{\mathcal{L}_\mathbb{N}} P(M', \cdot) P'(M, dM'), \quad M \in \mathcal{L}_\mathbb{N}. \tag{42}$$

By this operation, $\mathfrak{S}_\mathcal{L}$ is a semigroup with identity element given by $P(M, \cdot) = \mathbf{1}\{M \in \cdot\}$, the point mass at $M$, for every $M \in \mathcal{L}_\mathbb{N}$. Since $\mathfrak{S}_\mathcal{L}$ consists of exchangeable transition probability measures and $\mathcal{E}_\mathcal{L}$ consists of exchangeable probability measures, both on $\mathcal{L}_\mathbb{N}$, any $P \in \mathfrak{S}_\mathcal{L}$ acts on $\mu \in \mathcal{E}_\mathcal{L}$ by

$$(43) \quad (P \circ \mu)(\cdot) = (P\mu)(\cdot) := \int_{\mathcal{L}_\mathbb{N}} P(M, \cdot) \mu(dM), \quad M \in \mathcal{L}_\mathbb{N}.$$

**Definition 7.1** (Semigroup process). *Given a combinatorial Lévy process* $\mathbf{X} = (X_t, t \geq 0)$ *on* $\mathcal{L}_\mathbb{N}$, *we define its associated* semigroup process $\mathbf{P} = (P_t, t \geq 0)$ *on* $\mathfrak{S}_\mathcal{L}$ *by*

$$P_t(M, \cdot) = \|X_t\|(\{\Delta \in \mathcal{L}_\mathbb{N} : M \triangle \Delta \in \cdot\}), \quad M \in \mathcal{L}_\mathbb{N}, \quad t \geq 0,$$

*provided* $\|X_t\|$ *exists for all* $t \geq 0$.

Almost sure existence of the semigroup process $\mathbf{P} = (P_m, m \geq 0)$ associated to any discrete time exchangeable combinatorial Lévy process $\mathbf{X}$ is guaranteed by almost sure existence of $\|\mathbf{X}\|$. We define the increments of $\mathbf{X}$ to be $\Delta_m = X_m \triangle X_{m-1}$ for every $m = 1, 2, \ldots$. Each increment is exchangeable and, thus, $\|\Delta_m\|$ exists for every $m \geq 1$ and determines an exchangeable transition probability as in (41). We also see that $\|\Delta_{m+1} \triangle \Delta_m\| = \|\Delta_{m+1}\| \circ \|\Delta_m\|$ for all $m \geq 1$, where the righthand side is regarded as a composition of transition probabilities as in (42).

To see this explicitly, suppose $\|\Delta_{m+1}\| = W'$ and $\|\Delta_m\| = W$, then

$$\|\Delta_{m+1}\| \circ \|\Delta_m\|(\cdot) = \int_{\mathcal{L}_\mathbb{N}} W'(\{M' \in \mathcal{L}_\mathbb{N} : M \triangle M' \in \cdot\}) W(dM).$$

The measures $\|\Delta_m\|$, $\|\Delta_{m+1}\|$, and $\|\Delta_{m+1} \triangle \Delta_m\|$ are determined by how they behave on the $\pi$-system of events of form

$$\{M \in \mathcal{L}_\mathbb{N} : M|_{[n]} = A\}$$

for $A \in \mathcal{L}_{[n]}$, $n \in \mathbb{N}$. In this case, we let $A \in \mathcal{L}_{[k]}$ and note that

$$\begin{aligned}
\|\Delta_{m+1} \triangle \Delta_m\|(\{M \in \mathcal{L}_\mathbb{N} : M|_{[k]} = A\}) &= \delta(A; \Delta_{m+1} \triangle \Delta_m) \\
&= \lim_{n \to \infty} \frac{1}{n^{\downarrow k}} \sum_{\varphi:[k] \to [n]} \mathbf{1}\{(\Delta_{m+1} \triangle \Delta_m)^\varphi = A\} \\
&= \lim_{n \to \infty} \frac{1}{n^{\downarrow k}} \sum_{\varphi:[k] \to [n]} \mathbf{1}\{\Delta_{m+1}^\varphi \triangle \Delta_m^\varphi = A\} \\
&= \sum_{B \in \mathcal{L}_{[k]}} \lim_{n \to \infty} \frac{1}{n^{\downarrow k}} \sum_{\varphi:[k] \to [n]} \mathbf{1}\{\Delta_m^\varphi = B\} \mathbf{1}\{\Delta_{m+1}^\varphi = A \triangle B\}.
\end{aligned}$$

We now define

$$Z_n(B) := \frac{1}{n^{\downarrow k}} \sum_{\varphi:[k] \to [n]} \mathbf{1}\{\Delta_m^\varphi = B\} \mathbf{1}\{\Delta_{m+1}^\varphi = A \triangle B\}, \quad n \in \mathbb{N},$$

along with the $\sigma$-field $\mathcal{F}_n = \sigma\langle Z_{n+1}, Z_{n+2}, \ldots \rangle$ for each $n \in \mathbb{N}$. For any injection $\varphi : [k] \to [n]$, exchangeability and independence of the increments $\Delta_m$ and $\Delta_{m+1}$ implies

$$\begin{aligned}
\mathbb{P}\{\Delta_m^\varphi = B \text{ and } \Delta_{m+1}^\varphi = A \triangle B \mid \mathcal{F}_n\} &= \mathbb{P}\{\Delta_m|_{[k]} = B \text{ and } \Delta_{m+1}|_{[k]} = A \triangle B \mid \mathcal{F}_n\} \\
&= \mathbb{P}\{\Delta_m|_{[k]} = B \mid \mathcal{F}_n\} \mathbb{P}\{\Delta_{m+1}|_{[k]} = A \triangle B \mid \mathcal{F}_n\}.
\end{aligned}$$

Therefore,

$$\begin{aligned}
Z_{n+1}(B) &= \mathbb{E}\left(\frac{1}{(n+1)^{\downarrow k}} \sum_{\psi:[k]\to[n+1]} \mathbf{1}\{\Delta_m^\psi = B\}\mathbf{1}\{\Delta_{m+1}^\psi = A \triangle B\} \mid \mathcal{F}_n\right) \\
&= \frac{1}{(n+1)^{\downarrow k}} \sum_{\psi:[k]\to[n+1]} \mathbb{E}(\mathbf{1}\{\Delta_m^\psi = B\} \mid \mathcal{F}_n)\mathbb{E}(\mathbf{1}\{\Delta_{m+1}^\psi = A \triangle B\} \mid \mathcal{F}_n) \\
&= \mathbb{P}\{\Delta_m|_{[k]} = B \mid \mathcal{F}_n\}\mathbb{P}\{\Delta_{m+1}|_{[k]} = A \triangle B \mid \mathcal{F}_n\}.
\end{aligned}$$

It follows that $E(Z_n(B) \mid \mathcal{F}_n) = Z_{n+1}(B)$ and $(Z_n(B), n \in \mathbb{N})$ is a reverse martingale with respect to the filtration $(\mathcal{F}_n, n \in \mathbb{N})$. By the reverse martingale convergence theorem, there exists a random variable $Z_\infty(B)$ such that $Z_n(B) \to Z_\infty(B)$ almost surely.

By the analogous argument, we have

$$\lim_{n\to\infty} \mathbb{P}\{\Delta_m|_{[k]} = B \mid \mathcal{F}_n\} \to \mathbb{P}\{\Delta_m|_{[k]} = B \mid \mathcal{F}_\infty\} = \delta(B; \Delta_m) \quad \text{and}$$

$$\lim_{n\to\infty} \mathbb{P}\{\Delta_{m+1}|_{[k]} = A \triangle B \mid \mathcal{F}_n\} \to \mathbb{P}\{\Delta_{m+1}|_{[k]} = A \triangle B \mid \mathcal{F}_\infty\} = \delta(A \triangle B; \Delta_{m+1})$$

It follows that

$$\begin{aligned}
\|\Delta_{m+1} \triangle \Delta_m\|(A) &= \delta(A; \Delta_{m+1} \triangle \Delta_m) \\
&= \sum_{B \in \mathcal{L}_{[k]}} \delta(B; \Delta_m)\delta(A \triangle B; \Delta_{m+1}) \\
&= \|\Delta_{m+1}\| \circ \|\Delta_m\|(A) \quad \text{a.s.}
\end{aligned}$$

The stationary and independent increments properties, therefore, imply that $P_{m+n} =_{\mathcal{D}} P'_m \circ P_n$ for all $m, n \geq 0$, where $\mathbf{P}' = (P'_m, m \geq 0)$ is an independent copy of $\mathbf{P}$. In particular, $\|\mathbf{X}\| = (\|X_m\|, m \geq 0)$ can be constructed by taking $Q_1, Q_2, \ldots$ to be i.i.d. exchangeable transition measures with the same distribution as $P_1$ in the semigroup process associated to $\mathbf{X}$ and putting $\|X_m\| = Q_m \circ \|X_{m-1}\|$ for each $m \geq 1$, where the action of $Q_m$ on $\|X_{m-1}\|$ is as defined in (43).

$\square$

*Proof of Theorem 4.15: Finite case.* By the same argument as in the beginning of the proof of the infinite case, the increments $\Delta_m$ must satisfy (40) and, therefore, must be governed by an exchangeable probability measure on $\mathcal{L}_{[n]}$. The characterization of $\mu$ by some unique $\mathbf{p} = (p_Y)_{Y \in \mathcal{UL}_{[n]}}$ follows from Proposition 3.12.

Now consider the projection $\langle \mathbf{X} \rangle_\cong = (\langle X_m \rangle_\cong, m \geq 0)$ into $\mathcal{UL}_{[n]}$. By standard conditions under which a function of a Markov chain is a Markov chain, see [10], $\langle \mathbf{X} \rangle_\cong$ is a Markov chain just in case the transition probabilities of $\mathbf{X}$ satisfy

$$\mathbb{P}\{\langle X_{m+1}\rangle_\cong = Y' \mid X_m = M\} = \mathbb{P}\{\langle X_{m+1}\rangle_\cong = Y' \mid X_m = M'\}$$

for all $M, M'$ such that $\langle M \rangle_\cong = \langle M' \rangle_\cong$. Now, suppose $\langle X_m \rangle_\cong = Y$ and take $M, M' \in \mathcal{L}_{[n]}$ such that $\langle M \rangle_\cong = \langle M' \rangle_\cong = Y$. Then

$$\begin{aligned}
\mathbb{P}\{\langle X_{m+1}\rangle_\cong = Y' \mid X_m = M\} &= \sum_{M^* \in \mathcal{L}_{[n]}: \langle M^*\rangle_\cong = Y'} \mathbb{P}\{X_{m+1} = M^* \mid X_m = M\} \\
&= \sum_{M^* \in \mathcal{L}_{[n]}: \langle M^*\rangle_\cong = Y'} \mathbb{P}\{X_{m+1} = M^{*\sigma} \mid X_m = M^\sigma\}
\end{aligned}$$

for all permutations $\sigma : [n] \to [n]$ by exchangeability of $\mathbf{X}$. Let $\sigma : [n] \to [n]$ be any permutation such that $M^\sigma = M'$. (There must be at least one by the definition of $\langle M \rangle_\cong$ in (19) and our assumption that $\langle M \rangle_\cong = \langle M' \rangle_\cong$.) Then

$$
\begin{aligned}
\sum_{M^* \in \mathcal{L}_{[n]} : \langle M^* \rangle_\cong = Y'} \mathbb{P}\{X_{m+1} = M^{*\sigma} \mid X_m = M^\sigma\} &= \sum_{M^* \in \mathcal{L}_{[n]} : \langle M^* \rangle_\cong = Y'} \mathbb{P}\{X_{m+1} = M^{*\sigma} \mid X_m = M'\} \\
&= \sum_{M^* \in \mathcal{L}_{[n]} : \langle M^{*\sigma^{-1}} \rangle_\cong = Y'} \mathbb{P}\{X_{m+1} = M^* \mid X_m = M'\} \\
&= \sum_{M^* \in \mathcal{L}_{[n]} : \langle M^* \rangle_\cong = Y'} \mathbb{P}\{X_{m+1} = M^* \mid X_m = M'\} \\
&= \mathbb{P}\{\langle X_{m+1} \rangle_\cong = Y' \mid X_m = M'\}.
\end{aligned}
$$

It follows that $\langle \mathbf{X} \rangle_\cong = (\langle X_m \rangle_\cong, m \geq 0)$ is a Markov chain on $\mathcal{UL}_{[n]}$ with transition probabilities defined by

$$\mathbb{P}\{\langle X_{m+1} \rangle_\cong = Y' \mid \langle X_m \rangle_\cong = Y\} = \sum_{M' \in \mathcal{L}_{[n]} : \langle M' \rangle_\cong = Y'} \mathbb{P}\{X_{m+1} = M' \mid X_m = M\}$$

for any $M \in \mathcal{L}_{[n]}$ such that $\langle M \rangle_\cong = Y$.

□

□

7.2. **Continuous time combinatorial Lévy processes.** We first establish the Lévy property for the finite restrictions of any combinatorial Lévy process.

*Proof of Proposition 4.9.* The finite restrictions of any combinatorial Lévy process must also have stationary, independent increments and càdlàg sample paths, by the usual characterization of stochastic processes through their finite restrictions and the definition of the increment operator. Conversely, suppose $\mathbf{X} = (X_t, t \geq 0)$ is a stochastic process on $\mathcal{L}_\mathbb{N}$ whose finite restrictions are finite state space Lévy processes. Then $\mathbf{X}$ has càdlàg paths by the definition of the product discrete topology. Moreover, the increments of $\mathbf{X}$ are determined by the sequence of finite state space increments, so that the stationary and independent increments properties must also hold for $\mathbf{X}$. □

*Proof of Corollary 4.10.* By Proposition 4.9, each of the finite state space sample paths of $\mathbf{X}$ is also a Lévy process. Thus, each $\mathbf{X}^{[n]} = (X_t|_{[n]}, t \geq 0)$ has stationary, independent increments and càdlàg sample paths. For every $n \in \mathbb{N}$, $\mathcal{L}_{[n]}$ is a finite state space, and the càdlàg paths assumption implies that $\mathbf{X}^{[n]}$ has strictly positive hold times in all states it visits. By the Stone–Weierstrass theorem for compact Hausdorff spaces, for example, [20, Theorem 2.4.11],

$$C = \{g : \mathcal{L}_\mathbb{N} \to \mathbb{R} : \text{there exists } n \geq 1 \text{ such that } M|_{[n]} = M'|_{[n]} \text{ implies } g(M) = g(M')\}$$

is dense in the space of continuous, bounded functions $\mathcal{L}_\mathbb{N} \to \mathbb{R}$.

By Proposition 4.9, the Feller property is easily verified for each finite restriction $\mathbf{X}^{[n]}$ on account of the càdlàg paths assumption and the fact that $\mathcal{L}_{[n]}$ has a discrete topology. The Feller property of $\mathbf{X}$ now follows readily since the Feller property for each $\mathbf{X}^{[n]}$ implies that the conditions are satisfied by $\mathbf{X}$ for all functions in the dense set $C$ above.

□

Let $\mathbf{X} = (X_t, t \geq 0)$ be a combinatorial Lévy process on $\mathcal{L}_\mathbb{N}$. By Corollary 4.10, $\mathbf{X}$ has the Feller property and, therefore, its transition law is determined by the infinitesimal jump rates

$$Q(M, dM') := \lim_{t \downarrow 0} \frac{1}{t} \mathbb{P}\{X_t \in dM' \mid X_0 = M\} \quad M \neq M'.$$

By the stationary increments property, the jump rate from $M$ into $dM'$ depends only on the increment $\triangle(M, M')$. Thus, we can define a measure

(44) $$\mu(d\triangle) := \begin{cases} Q(\mathbf{0}_\mathbb{N}^\mathcal{L}, d\triangle), & \triangle \neq \mathbf{0}_\mathbb{N}^\mathcal{L}, \\ 0, & \triangle = \mathbf{0}_\mathbb{N}^\mathcal{L}. \end{cases}$$

**Proposition 7.2.** *The measure $\mu$ defined in* (44) *satisfies* (22).

*Proof.* The lefthand side of (22) is plain by (44), which requires $\mu(\mathbf{0}_\mathbb{N}^\mathcal{L}) = 0$. The righthand side follows from Proposition 4.9 as we now show. By construction, $\mathbf{X}^{[n]}$ has infinitesimal jump rates

$$Q_n(M, M') := Q(M^*, \{M'' \in \mathcal{L}_\mathbb{N} : M''|_{[n]} = M'\}), \quad M' \neq M \in \mathcal{L}_{[n]},$$

for any $M^* \in \{M'' \in \mathcal{L}_\mathbb{N} : M''|_{[n]} = M\}$, for each $n \in \mathbb{N}$. We interpret $Q_n(M, M')$ as the rate at which $\mathbf{X}^{[n]}$ jumps from $M$ to $M'$, and so we put $Q_n(M, M) = 0$ for all $M \in \mathcal{L}_{[n]}$. Since $\mathcal{L}_{[n]}$ is finite and $Q_n(M, M') < \infty$ for all $M' \neq M$, we have

$$\infty > Q_n(M, \mathcal{L}_{[n]}) = \mu(\{M^* \in \mathcal{L}_\mathbb{N} : M^*|_{[n]} \neq \mathbf{0}_{[n]}^\mathcal{L}\})$$

and the righthand side of (22) holds. □

*Proof of Theorem 4.13.* Let $\mathbf{X} = (X_t, t \geq 0)$ be a combinatorial Lévy process with infinitesimal jump measure $\mu$ defined in (44). Let $\mathbf{X}_\mu^* = (X_t^*, t \geq 0)$ be a $\mu$-canonical Lévy process, as in Section 4.1. Since $\mathbf{X}_\mu^*$ is constructed from the finite state space processes $\mathbf{X}_\mu^{*[n]}$, its jump rates are determined by

(45) $$\mu^{(n)}(\triangle) = \begin{cases} \mu(\{M \in \mathcal{L}_\mathbb{N} : M|_{[n]} = \triangle\}), & \triangle \neq \mathbf{0}_{[n]}^\mathcal{L}, \\ 0, & \triangle = \mathbf{0}_{[n]}^\mathcal{L}. \end{cases}$$

We define $\mu^*$ on events of the form $\{M^* \in \mathcal{L}_\mathbb{N} : M^*|_{[n]} = M\}$, for $M \in \mathcal{L}_{[n]} \setminus \{\mathbf{0}_{[n]}^\mathcal{L}\}$, for every $n \in \mathbb{N}$, by

$$\mu^*(\{M^* \in \mathcal{L}_\mathbb{N} : M^*|_{[n]} = M\}) = \mu^{(n)}(M).$$

Such events comprise a generating $\pi$-system of the Borel $\sigma$-field on $\mathcal{L}_\mathbb{N} \setminus \{\mathbf{0}_\mathbb{N}^\mathcal{L}\}$ and $\mu^*$ is additive by construction. Carathéodory's extension theorem implies a unique extension of $\mu^*$ to $\mathcal{L}_\mathbb{N} \setminus \{\mathbf{0}_\mathbb{N}^\mathcal{L}\}$. Putting $\mu^*(\{\mathbf{0}_\mathbb{N}^\mathcal{L}\}) = 0$ gives a unique measure on $\mathcal{L}_\mathbb{N}$. Since $\mu^*$ coincides with $\mu$ on the generating $\pi$-system, we must have $\mu^* = \mu$. This completes the proof. □

We can now immediately deduce the Lévy–Itô characterization for infinite exchangeable combinatorial Lévy processes, which we restate for the reader's convenience.

**Theorem 4.18.** *Let $\mathcal{L} = (i_1, \ldots, i_k)$ be any signature with $i_k \geq 1$ and $\mathbf{X} = (X_t, t \geq 0)$ be an exchangeable combinatorial Lévy process on $\mathcal{L}_\mathbb{N}$. Then there exists a unique measure $\nu_0$ on $\mathcal{E}_\mathcal{L}$ satisfying*

(46) $$\nu_0(\{\mathbf{0}_\mathcal{L}\}) = 0 \quad \text{and} \quad \int_{\mathcal{E}_\mathcal{L}} (1 - W(\{\mathbf{0}_{[i_k]}^\mathcal{L}\}))\nu_0(dW) < \infty,$$

*and measures $\mu_\lambda$ on $\mathcal{L}_\mathbb{N}$ for every $\lambda \vdash q$, $q = 1, \ldots, i_k$, satisfying* (30) *and* (31) *such that*

$$\mu = \nu_0^* + \sum_{q=1,\ldots,i_k} \sum_{\lambda \vdash q} \mu_\lambda^*, \tag{47}$$

*where $\mu_\lambda^*$ is defined in* (32).

*Proof of Theorem 4.18.* The proof follows from Theorems 4.6, 4.13, and 6.3. □

7.3. **Proof of Theorem 4.19.** Theorem 4.19 characterizes the existence and behavior of the limiting process $\|\mathbf{X}\|$. The theorem has two key parts. We first show that the projection of an exchangeable Lévy process $\mathbf{X} = (X_t, t \geq 0)$ on $\mathcal{L}_\mathbb{N}$ exists almost surely at all time points. We then show that $\|\mathbf{X}\|$ is itself a Feller process whose discontinuities we characterize in Theorem 4.20 and Corollary 4.21 in terms of the jumps of the covering process $\mathbf{X}$.

7.3.1. *Existence.* By exchangeability of $\mathbf{X}$, Proposition 6.2 implies that $\|X_t\|$ exists for any countable collection of times $t$. This argument does not generalize to existence at all times $t$, as there are uncountably many of them. To prove existence of $\|\mathbf{X}\|$ simultaneously at all times, we show that it exists at all $t \in [0,1]$ with probability 1. We deduce existence for all $t \in [0,\infty)$ by partitioning $[0,\infty) = \bigcup_{n \geq 1}[n-1, n)$, proving existence of $\|\mathbf{X}\|$ for all $t \in [n-1, n)$ for each $n \geq 1$ by time homogeneity of combinatorial Lévy processes, and deducing existence of $\|\mathbf{X}\|$ for all $t \in [0,\infty)$ by countable additivity of probability measures.

For every $m \in \mathbb{N}$, we define the upper and lower homomorphism densities of $A \in \mathcal{L}_{[m]}$ in $M \in \mathcal{L}_\mathbb{N}$ by

$$\delta^+(A; M) := \limsup_{n \to \infty} \frac{1}{n^{\downarrow m}} \sum_{\varphi:[m] \to [n]} \mathbf{1}\{M|_{[n]}^\varphi = A\} \quad \text{and}$$

$$\delta^-(A; M) := \liminf_{n \to \infty} \frac{1}{n^{\downarrow m}} \sum_{\varphi:[m] \to [n]} \mathbf{1}\{M|_{[n]}^\varphi = A\}.$$

We define the upper and lower combinatorial limits, respectively, by

$$\|M\|^+ := (\delta^+(A; M), A \in \bigcup_{m \in \mathbb{N}} \mathcal{L}_{[m]}) \quad \text{and}$$

$$\|M\|^- := (\delta^-(A; M), A \in \bigcup_{m \in \mathbb{N}} \mathcal{L}_{[m]}).$$

Since the limits inferior and superior always exist, the upper and lower limits of $M$ are well defined. The limit $\|M\|$ exists if and only if $\|M\|^+ = \|M\|^-$, so we must show that $\|X\|^+ = (\|X_t\|^+, 0 \leq t \leq 1)$ and $\|X\|^- = (\|X_t\|^-, 0 \leq t \leq 1)$ coincide with probability 1.

By the canonical construction of $\mathbf{X}$ from a time homogeneous Poisson point process $\Delta^*$ with intensity $dt \otimes \mu$, there is probability 0 of a discontinuity at any given time $t \in [0,1]$. By the càdlàg paths assumption, each finite restriction $\mathbf{X}^{[n]}$ of $\mathbf{X}$ experiences at most finitely many discontinuities in $[0,1]$. Moreover, by the definition of the increments operator, $\mathbf{X}^{[n]}$ experiences a discontinuity in any interval $[s,t]$, $s < t$, if and only if there is an atom time $u \in [s,t]$ in $\Delta^*$ such that $\Delta_u|_{[n]} \neq \mathbf{0}_{[n]}^\mathcal{L}$. Thus, for any $s < t$,

$$\mathbb{P}\{\mathbf{X}^{[n]} \text{ is discontinuous on } [s,t]\} = 1 - \exp\{-(t-s)\mu^{(n)}(\mathcal{L}_{[n]})\}$$

for $\mu^{(n)}$ as defined in (45). In particular, for any $\varepsilon > 0$,

$$\mathbb{P}\{\mathbf{X}^{[n]} \text{ is discontinuous on } [s,t]\} < \varepsilon$$

as long as $t - s < -\log(1-\varepsilon)/\mu^{(n)}(\mathcal{L}_{[n]})$.

Thus, for every $m \in \mathbb{N}$ and every $\varepsilon > 0$, we define $s = -\frac{1}{2}\log(1-\varepsilon)/\mu^{(m)}(\mathcal{L}_{[m]})$ and partition $[0, 1]$ into finitely many non-overlapping subintervals $[0, s), [s, 2s), \ldots, [\lfloor 1/s \rfloor s, 1]$ such that

$$\mathbb{P}\{\Delta^* \text{ has an atom } (t, \Delta_t) \text{ in } [ks, (k+1)s] \text{ for which } \Delta_t|_{[m]} \neq \mathbf{0}^{\mathcal{L}}_{[m]}\} < \varepsilon \quad \text{for every } k = 0, 1, \ldots, \lfloor 1/s \rfloor.$$

Since the action of relabeling is ergodic for exchangeable processes, the law of large numbers implies

$$\lim_{n \to \infty} \frac{1}{n^{\downarrow m}} \sum_{\varphi:[m] \to [n]} \mathbf{1}\{(X_t|^\varphi_{[n]}, 0 \leq t \leq 1) \text{ is discontinuous on } [ks, (k+1)s]\} < \varepsilon,$$

for every $k = 0, 1, \ldots, \lfloor 1/s \rfloor$. Thus, for each $k = 0, 1, \ldots, \lfloor 1/s \rfloor$, the upper and lower homomorphism densities of any $A \in \mathcal{L}_{[m]}$ in $X_t$ cannot vary by more than $\varepsilon$ over $[ks, (k+1)s)$. Furthermore, the upper and lower densities are equal on the fixed set of endpoints of $0, s, 2s, \ldots, \lfloor 1/s \rfloor s, 1$, implying

$$\sup_{t \in [0,1]} |\delta^+(A; X_t) - \delta^-(A; X_t)| \leq 2\varepsilon \quad \text{a.s.},$$

for every $A \in \mathcal{L}_{[m]}$, for every $m \in \mathbb{N}$. Since $\varepsilon > 0$ is arbitrary, it follows that $\delta^+(A; X_t) = \delta^-(A; X_t)$ simultaneously for all $t \in [0, 1]$ with probability 1. Countable additivity implies that $\|X_t\|^+ = \|X_t\|^-$ simultaneously for all $t \in [0, 1]$ with probability 1. By stationarity, the same argument applies to establish that $\|X_t\|^+ = \|X_t\|^-$ simultaneously for all $t \in [k, k+1]$ with probability 1 for every $k = 1, 2, \ldots$. Countable additivity implies $\|X_t\|^+ = \|X_t\|^-$ simultaneously for all $t \geq 0$ with probability 1.

7.3.2. *The Feller property.* To prove the Feller property for $\|\mathbf{X}\|$, we recall the metric

$$d(W, W') := \sum_{n \in \mathbb{N}} 2^{-n} \sum_{A \in \mathcal{L}_{[n]}} |W(A) - W'(A)|, \quad W, W' \in \mathcal{E}_\mathcal{L}.$$

The semigroup operation in (43) is Lipschitz continuous in this metric.

Almost sure existence of the semigroup process $\mathbf{P} = (P_t, t \geq 0)$ associated to $\mathbf{X}$ follows from almost sure existence of $\|\mathbf{X}\|$ for every exchangeable combinatorial Lévy process $\mathbf{X}$. By the analogous argument to the discrete time case, the stationary and independent increments properties of $\mathbf{X}$ imply that $\mathbf{P}$ satisfies $P_{t+s} =_\mathcal{D} P'_t \circ P_s$, for all $s, t \geq 0$, where $\mathbf{P}' = (P'_t, t \geq 0)$ is an independent, identically distributed copy of $\mathbf{P}$.

By the Poisson point process construction of $\mathbf{X}$, we can couple $\mathbf{X}$, $\|\mathbf{X}\|$, and $\mathbf{P}$ so that $\|X_t\| = P_t \circ \|X_0\|$ for all $t \geq 0$. Let $\mathbf{Q} = (Q_t)_{t \geq 0}$ be the semigroup of $\|\mathbf{X}\|$, that is,

$$Q_t g(W) = \mathbb{E}(g(\|X_t\|) \mid \|X_0\| = W).$$

We must establish the two conditions:

(i) $W \mapsto Q_t g(W)$ is continuous for all $t > 0$ and
(ii) $\lim_{t \downarrow 0} Q_t g(W) = g(W)$ for all $W \in \mathcal{E}_\mathcal{L}$.

Part (i) follows by continuity of $g$ and Lipschitz continuity of the semigroup action. Specifically, fix $W \in \mathcal{E}_\mathcal{L}$ and let $(W_n, n \geq 1)$ be a sequence in $\mathcal{E}_\mathcal{L}$ such that $W_n \to W$. For any bounded, continuous $g : \mathcal{E}_\mathcal{L} \to \mathbb{R}$, we have

$$\begin{aligned}
\lim_{n\to\infty} E(g(\|X_t\|) \mid \|X_0\| = W_n) &= \lim_{n\to\infty} E(g(P_t \circ \|X_0\|) \mid \|X_0\| = W_n) \\
&= \lim_{n\to\infty} E(g(P_t \circ W_n)) \\
&= E(g(P_t \circ W)),
\end{aligned}$$

because the composition of the bounded, continuous functions $g$ and $P_t$ is again bounded and continuous.

For (ii), we must have

(48) $$\mathbb{E}g(\|X_t\| \mid \|X_0\| = W) \to g(W) \quad \text{as } t \downarrow 0.$$

By the Feller property for $\mathbf{X}$, we know that $X_0 \triangle \Delta_t \to_\mathcal{D} X_0$ as $t \downarrow 0$. We deduce (48) from the implication (i)$\Rightarrow$(ii) in the following lemma.

**Lemma 7.3.** *For each $n \in \mathbb{N}$, let $X^{(n)}$ be an infinite exchangeable random $\mathcal{L}$-structure and let $\|X^{(n)}\|$ be its associated combinatorial limit. The following are equivalent.*

(i) $X^{(n)} \to_\mathcal{D} X^{(\infty)}$ as $n \to \infty$.
(ii) $\|X^{(n)}\| \to_\mathcal{D} \|X^{(\infty)}\|$ as $n \to \infty$.

*Proof of (ii)$\Longrightarrow$(i).* Suppose that $\|X^{(n)}\| \to_\mathcal{D} \|X^{(\infty)}\|$. Since $\mathcal{E}_\mathcal{L}$ is a compact metric space under (18), Skorokhod's representation theorem allows us to put the sequence $(\|X^{(n)}\|, n \in \mathbb{N})$ on the same probability space such that $\|X^{(n)}\| \to \|X^{(\infty)}\|$ a.s. By definition, $\|X^{(n)}\| = (\delta(A; X^{(n)}), A \in \bigcup_{m\in\mathbb{N}} \mathcal{L}_{[m]})$ for each $n \in \mathbb{N}$, and so $\|X^{(n)}\| \to \|X^{(\infty)}\|$ a.s. implies $\delta(A; X^{(n)}) \to \delta(A; X^{(\infty)})$ a.s. for every $A \in \bigcup_{m\in\mathbb{N}} \mathcal{L}_{[m]}$.

For each $m \in \mathbb{N}$, fix an ordering of $\mathcal{L}_{[m]}$, say, $A_{m,1}, A_2, \ldots, A_{m,|\mathcal{L}_{[m]}|}$, and define $\delta_{m,\bullet}(A_{0,j}; X^{(n)}) = 0$ and

$$\delta_{m,\bullet}(A_{m,j}; X^{(n)}) := \delta(A_{m,1}; X^{(n)}) + \delta(A_{m,2}; X^{(n)}) + \cdots + \delta(A_{m,j}; X^{(n)})$$

for each $j = 1, \ldots, |\mathcal{L}_{[m]}|$, so that $\delta_{m,\bullet}(A_j; X^{(n)}) \to \delta_{m,\bullet}(A_j; X^{(\infty)})$ a.s. for all $A_j \in \mathcal{L}_{[m]}$, for all $m \in \mathbb{N}$. For each $M \in \mathcal{L}_{[m]}$ and $k \geq m$, we define the conditional homomorphism density

$$\delta(A; X^{(n)} \mid M) = \begin{cases} \delta(A; X^{(n)})/\delta(M; X^{(n)}), & A|_{[m]} = M \text{ and } \delta(M; X^{(n)}) > 0, \\ 0, & \text{otherwise} \end{cases}, \quad A \in \mathcal{L}_{[k]}.$$

We then put

$$\delta_{m,\bullet}(A_{m,j}; X^{(n)} \mid M) = \delta(A_{m,1}; X^{(n)} \mid M) + \cdots \delta(A_{m,j}; X^{(n)} \mid M).$$

Now let $\xi_1, \xi_2, \ldots$ be i.i.d. Uniform$[0, 1]$ random variables. For every $n \in \mathbb{N}$, we construct $Y^{(n)}$ from the same distribution as $X^{(n)}$ by putting $Y^{(n)}|_{[1]} = A_j$ if $\delta_{1,\bullet}(A_{1,j-1}; X^{(n)}) \leq \xi_1 < \delta_{1,\bullet}(A_{1,j}; X^{(n)})$, with the convention that $\delta_{m,\bullet}(A_{m,0}; X^{(n)}) \equiv 0$ for all $m$. Given $Y^{(n)}|_{[m]} = M$, we put $Y^{(n)}|_{[m+1]} = A_{m+1,j} \in \mathcal{L}_{[m+1]}$ if $\delta_{m+1,\bullet}(A_{m+1,j-1}; X^{(n)} \mid M) \leq \xi_{m+1} < \delta_{m+1,\bullet}(A_{m+1,j}; X^{(n)} \mid M)$. Proceeding in this way produces a compatible sequence $(Y^{(n)}|_{[m]}, m \geq 1)$ of finite structures and, thus, a random structure $Y^{(n)} =_\mathcal{D} X^{(n)}$.

We construct each $Y^{(n)}, n \in \mathbb{N}$, from the same i.i.d. sequence of uniform random variables $\xi_1, \xi_2, \ldots$. Since we have assumed $\|X^{(n)}\| \to \|X^{(\infty)}\|$ a.s., we must have $Y^{(n)} \to Y^{(\infty)}$ a.s., where $Y^{(\infty)} =_\mathcal{D} X^{(\infty)}$. It follows that $X^{(n)} \to_\mathcal{D} X^{(\infty)}$.

$\square$

*Proof of (i)⟹(ii).* Suppose that $X^{(n)} \to_{\mathcal{D}} X^{(\infty)}$. Under our metric $d(\cdot, \cdot)$ from (18), $\mathcal{E}_{\mathcal{L}}$ is a Polish space and, thus, the space of probability measures on $\mathcal{E}_{\mathcal{L}}$ is compact by Prokhorov's theorem; see, e.g., [20, Chapter 11]. It follows that we can extract a subsequence $\|X^{(n_k)}\|$ such that $\|X^{(n_k)}\| \to_w \|\tilde{X}^{(\infty)}\|$ as $k \to \infty$, where $\to_w$ denotes weak convergence. By assumption, $X^{(n)} \to_{\mathcal{D}} X^{(\infty)}$ and each $X^{(n)}$ is conditionally distributed according to $\|X^{(n)}\|$; thus, the subsequence $X^{(n_k)}$ also converges to $X^{(\infty)}$, whose law must be given by $\|\tilde{X}^{(\infty)}\|$. Thus, the limit $\|\tilde{X}^{(\infty)}\|$ does not depend on the choice of subsequence and we must have $\|X^{(n)}\| \to \|\tilde{X}^{(\infty)}\| =_{\mathcal{D}} \|X^{(\infty)}\|$. □

### 7.4. **Proof of Theorem 4.20.**

7.4.1. *Classification of discontinuities.* By Theorem 4.18, the infinitesimal jump rates of any exchangeable combinatorial Lévy process $\mathbf{X} = (X_t, t \geq 0)$ are determined by an exchangeable measure

$$\mu = v_0^* + \sum_{q=1}^{i_k} \sum_{\lambda \vdash q} \mu_\lambda^*,$$

for $v_0$ satisfying (33) and $\mu_\lambda^*$ defined as (32) for $\mu_\lambda$ satisfying (22), (30), and (31) for each $\lambda \vdash q$, $q = 1, \ldots, i_k$. By Theorem 4.13, $\mathbf{X}$ can be constructed from a Poisson point process $\mathbf{\Delta} = \{(t, \Delta_t)\} \subset [0, \infty) \times \mathcal{L}_{\mathbb{N}}$ with intensity $dt \otimes \mu$. By the lefthand side of (33), $v_0^*$-almost every $M \in \mathcal{L}_{\mathbb{N}}$ has $\|M\| \neq \mathbf{0}_{\mathcal{L}}$, implying that the measures $v_0^*$ and $(\mu_\lambda^*, \lambda \vdash q = 1, \ldots, i_k)$ comprising $\mu$ are mutually singular. By properties of Poisson point processes, we can construct $\mathbf{\Delta}$ alternatively as a superposition of independent Poisson point processes $\mathbf{\Delta}_0$ and $\mathbf{\Delta}_\lambda$ for each $\lambda \vdash q = 1, \ldots, i_k$, all on $[0, \infty) \times \mathcal{L}_{\mathbb{N}}$, so that $\mathbf{\Delta}_0$ has intensity $dt \otimes v_0^*$ and $\mathbf{\Delta}_\lambda$ has intensity $dt \otimes \mu_\lambda^*$ for each $\lambda$.

In this way, we define $\mathbf{X}^0 = (X_t^0, t \geq 0)$ to be the process on $\mathcal{L}_{\mathbb{N}}$ constructed as in (4) from $\mathbf{\Delta}_0$ and, for each $\lambda$, $\mathbf{X}^\lambda = (X_t^\lambda, t \geq 0)$ to be the process on $\mathcal{L}_{\mathbb{N}}$ constructed as in (4) from the Poisson point process $\mathbf{\Delta}_\lambda$.

**Proposition 7.4.** *Let $\mathbf{X}_\mu^* = (X_t^*, t \geq 0)$ be a $\mu$-canonical Lévy process for some exchangeable measure $\mu = v_0^* + \sum_{q=1}^{i_k} \sum_{\lambda \vdash q} \mu_\lambda^*$ satisfying the conditions of Theorem 4.18. Let $\mathbf{X}^0$ and $\mathbf{X}^\lambda$ be as defined above. Then $\mathbf{X}_\mu^* =_{\mathcal{D}} \triangle_{0 \leq q \leq i_k, \lambda \vdash q} \mathbf{X}^\lambda$, where $\triangle_{0 \leq q \leq i_k, \lambda \vdash q} \mathbf{X}^\lambda$ is the process $(Y_t, t \geq 0)$ defined by*

$$Y_t := \triangle_{0 \leq q \leq i_k, \lambda \vdash q} X_t^\lambda, \quad t \geq 0. \tag{49}$$

*Proof.* By the superposition property of Poisson point processes, $\mathbf{\Delta}$ with intensity $dt \otimes \mu$ can be obtained as the superposition of $\mathbf{\Delta}_\lambda$ for $\lambda \vdash q = 0, 1, \ldots, i_k$. The $\triangle$ operator is commutative and associative, allowing the construction in (49). □

Now consider an exchangeable, dissociated combinatorial Lévy process $\mathbf{X} = (X_t, t \geq 0)$ with characteristic measure $\mu$. Since the process $\mathbf{X}$ is dissociated then $X_t$ is marginally dissociated for every $t \geq 0$. Let

$$Q(M, dM') := \mu(\{\Delta \in \mathcal{L}_{\mathbb{N}} : M \triangle \Delta \in dM'\}), \quad M, M' \in \mathcal{L}_{\mathbb{N}},$$

be the infinitesimal jump rates of $\mathbf{X}$ induced by $\mu$. The distribution of $X_t$ is determined by the dissociated infinitesimal rate measure $\mu$ so that

$$P_t(M, \cdot) = \mathbb{P}\{X_t \in \cdot \mid X_0 = M\} = e^{-tQ}(M, \cdot) = \sum_{k=0}^{\infty} \frac{t^k}{k!} Q^{(k)}(M, \cdot),$$

where $Q^{(k)}$ is the $k$-fold composition of $Q$ with itself. (For every $n \in \mathbb{N}$, we define $Q_n$ by

$$Q_n(A, A') = Q(M, \{M' \in \mathcal{L}_{\mathbb{N}} : M'|_{[n]} = A\}), \quad A, A' \in \mathcal{L}_{[n]},$$

for any $M \in \mathcal{L}_{\mathbb{N}}$ such that $M|_{[n]} = A$. Since $\mathcal{L}_{[n]}$ is finite, we can arrange $Q_n$ in a matrix with diagonal entries $Q_n(A, A) = -Q_n(A, \mathcal{L}_{[n]} \setminus \{A\})$. The $k$-fold measure $Q^{(k)}$ is determined by taking the $k$-fold matrix product $Q_n^k$ of $Q_n$ and putting

$$Q^{(k)}(M, \{M' \in \mathcal{L}_{\mathbb{N}} : M'|_{[n]} = A'\}) = Q_n^k(M|_{[n]}, A'), \quad A' \in \mathcal{L}_{[n]}.)$$

Since $\mathbf{X}$ is dissociated and assumed to have $X_0 = \mathbf{0}_{\mathbb{N}}^{\mathcal{L}}$, the combinatorial limit $\|X_t\|$ is deterministic at every fixed $t > 0$; thus,

$$\begin{aligned} \delta(A; X_t) &= \mathbb{E}\delta(A; X_t) \\ &= \mathbb{P}\{X_t|_{[n]} = A\} \\ &= e^{tQ}(\mathbf{0}_{[n]}^{\mathcal{L}}, A) \end{aligned}$$

with probability 1.

Now, let $\mathbb{Q}$ denote the rational numbers and for each $q \in \mathbb{Q}$ define $Y_q = e^{qQ}(\mathbf{0}_{\mathbb{N}}^{\mathcal{L}}, \cdot)$. Then $(\|X_q\|, q \in \mathbb{Q}) = (Y_q, q \in \mathbb{Q})$ a.s. since $\mathbb{Q}$ is countable. We extend $(Y_q, q \in \mathbb{Q})$ to a process on all of $[0, \infty)$ by putting

$$Y_t = \lim_{q \downarrow t} Y_q = \lim_{q \downarrow t} e^{qQ}(\mathbf{0}_{\mathbb{N}}^{\mathcal{L}}, \cdot) = e^{tQ}(\mathbf{0}_{\mathbb{N}}^{\mathcal{L}}, \cdot)$$

for every $t \geq 0$.

By definition, $\mathbf{Y} = (Y_t, t \geq 0)$ is the continuous, deterministic path $t \mapsto e^{tQ}(\mathbf{0}_{\mathbb{N}}^{\mathcal{L}}, \cdot)$ on $\mathcal{E}_{\mathcal{L}}$. By the assumption that $\mathbf{X}$ is dissociated, we have $\mathbb{P}\{\|X_t\| = Y_t\} = 1$ for all fixed $t \geq 0$, so that $\mathbf{Y}$ is a version of $\|\mathbf{X}\|$. Furthermore, $\|\mathbf{X}\| =_{\mathcal{D}} \mathbf{Y}$ since $\mathbb{P}\{\|X_{t_j}\| = Y_{t_j}, j = 1, \ldots, r\} = 1$ holds for all finite sets of times $0 \leq t_1 < \cdots < t_r < \infty$.

To establish Corollary 4.21, we once again appeal to the decomposition in Proposition 7.4. We showed previously, in the proof of Theorem 4.15, that $\|X \triangle X'\| = \|X\| \circ \|X'\|$ a.s. whenever $X, X'$ are independent, exchangeable random $\mathcal{L}$-structures. From this, we let $\mathbf{X}^\lambda$ be the $\lambda$-component of $\mathbf{X}_\mu^*$ from Proposition 7.4. By exchangeability, $\|\mathbf{X}^\lambda\|$ exists a.s. for every $\lambda \vdash q = 0, 1, \ldots, i_k$. Moreover, $\mathbf{X}^\lambda$ is dissociated for each $\lambda \vdash q = 1, \ldots, i_k$ by definition of the $\mu_\lambda$ measures; whence $\|\mathbf{X}^\lambda\| =_{\mathcal{D}} (e^{tQ}(\mathbf{0}_{\mathbb{N}}^{\mathcal{L}}, \cdot), t \geq 0)$ by Theorem 4.20 above. It follows that

$$\begin{aligned} \|\mathbf{X}\| &=_{\mathcal{D}} \|\triangle_{\lambda \vdash q=0,1\ldots,i_k} \mathbf{X}^\lambda\| \\ &= (\|\triangle_{\lambda \vdash q=1,\ldots,i_k} X_t^\lambda \triangle X_t^0\|, t \geq 0) \\ &= (e^{t \sum_{\lambda \vdash q=1,\ldots,i_k} Q^\lambda}(\mathbf{0}_{\mathbb{N}}^{\mathcal{L}}, \cdot) \circ \|X_t^0\|, t \geq 0) \quad \text{a.s.} \end{aligned}$$

If $\|X_t^0\|$ is continuous at $t > 0$, then so is $\|\mathbf{X}\|$.

**Remark 7.5.** *The above argument classifies the discontinuities of $\|\mathbf{X}\|$ as a subset of the jumps from the $\nu_0^*$ measure in the decomposition* (34). *An analogous classification in the special case of graph-valued processes is claimed in* [15, Theorem 5]; *however, the argument given in* [15] *only proves the corresponding statement for each fixed time $t > 0$, which does not immediately extend to the uncountable set of all times $t > 0$.*

## 8. Closing remarks

8.1. **Applications.** Stochastic process models for dynamic combinatorial structures have a place in DNA sequencing, dynamic network modeling, combinatorial search algorithms, and much more. They also have potential for modeling certain composite structures, as discussed in Section 5.3. With an array of applications in mind, we have developed the theory of combinatorial Lévy processes in a general setting and proven some basic theorems about their behavior.

The description of combinatorial Lévy processes as a process with stationary, independent increments, as we have defined it here, facilitates statistical inference of the jump measure, tests for exchangeability, and tests for stationarity. A particularly important aspects of this theory is the ability to analyze combinatorial Lévy processes that are not exchangeable and which evolve on finite state spaces, neither of which is possible in other studies of combinatorial stochastic processes, for example, in [9, 15, 35]. The significance of this is pronounced in statistical applications, where the assumption of infinite exchangeability often carries with it additional baggage that may constrain the available models in a prohibitive way.

As a demonstration of potential statistical questions that can be handled by combinatorial Lévy processes but not other theories of combinatorial processes, suppose $\mathbf{X} = (X_m)_{m=0,1,\ldots,T}$ is an observation of sequence of structures in $\mathcal{L}_{[n]}$ for some signature $\mathcal{L}$. From $\mathbf{X}$, we compute the empirical jump distribution $\hat{\mu}$ by

$$(50) \qquad \hat{\mu}(M) := \frac{1}{T} \sum_{t=1}^{T} \mathbf{1}_M(X_t \triangle X_{t-1}), \quad M \in \mathcal{L}_{[n]},$$

where $\mathbf{1}_M(\cdot)$ is the unit mass at $M$. From this empirical measure, we can estimate an empirical jump measure $\hat{\mu}^{\mathrm{ex}}$ by

$$(51) \qquad \hat{\mu}^{\mathrm{ex}}(M) := \frac{1}{|\{M' \in \mathcal{L}_{[n]} : \langle M' \rangle_\cong = \langle M \rangle_\cong\}|} \sum_{M' \in \mathcal{L}_{[n]} : \langle M' \rangle_\cong = \langle M \rangle_\cong} \hat{\mu}(M'), \quad M \in \mathcal{L}_{[n]},$$

the measure obtained by averaging over equivalence classes. The goodness of fit for the exchangeable model in (51) can be compared to (50) by a Pearson chi-square test statistic to test for exchangeability.

Alternatively, in the special case of graph-valued processes, we can test the Lévy process assumption of stationarity against the general class of exchangeable Feller chains studied in [15]. Those processes, however, are limited to the assumption of infinite populations and would require additional data to estimate the transition measure.

## Acknowledgment

The author's work is partially supported by NSF grants CNS-1523785 and CAREER DMS-1554092.

Rutgers University, Department of Statistics, 110 Frelinghuysen Road, Piscataway, NJ 08854.